\documentclass[final,3p,times]{elsarticle}

\usepackage[normalem]{ulem}

\usepackage{amssymb}
\usepackage{amsmath}
\usepackage{amsthm}
\usepackage{subeqnarray}
\usepackage{mathrsfs}
\usepackage{rotating}
\usepackage{enumitem}
\setlist[enumerate]{noitemsep, topsep=2pt}
\setlist[itemize]{noitemsep, topsep=2pt}
\usepackage{bm}             
\usepackage{setspace}

\usepackage{algorithmic}
\usepackage[linesnumbered,ruled]{algorithm2e}

\usepackage[breaklinks,bookmarks=false]{hyperref}
\hypersetup{colorlinks, linkcolor=blue, citecolor=blue,
  urlcolor=blue, pagecolor=blue}

\usepackage{cleveref}
\usepackage{amstext}

\biboptions{sort&compress}

\DeclareMathOperator*{\Div}{div}
\DeclareMathOperator*{\rot}{rot}
\DeclareMathOperator*{\curl}{curl}

\newtheorem{theorem}{Theorem}[section]

\theoremstyle{remark}
\newtheorem{remark}[theorem]{Remark}

\newcommand\bifont {\bm}

\newcommand\matN{\mathbb{N}}

\newcommand\krr{{\mathscr R}}

\newcommand\bkb{\bifont{b}}
\newcommand\bkn{\bifont{n}}
\newcommand\bkq{\bifont{q}}
\newcommand\bkr{\bifont{r}}

\newcommand\Krr[1]{{\krr({#1})}}
\newcommand\KrrO{{\krr}}

\newcommand\pd {\partial}

\newcommand\vp {\varphi}
\newcommand\vph {\varphi_h}
\newcommand\psih {\psi_h}

\newcommand\RR {\mathbb{R}}
\newcommand\IR {\RR}
\newcommand\R  {\RR}    
\newcommand\ve {\varepsilon}
\newcommand\Om {{\Omega}}
\newcommand\oO {\overline{\om}}

\newcommand\gom{{\Gamma}}

\newcommand\dd {\mathrm{d}}

\newcommand\dx {{\,\dd x}}

\newcommand\dS {{\,\dd S}}
\newcommand\ds {{\,\dd s}}

\newcommand\SShp {S_{h}^p}

\newcommand\Shp {\Wh}

\newcommand\Shpl {\Whm}

\renewcommand\u {u}

\newcommand\z {z}

\newcommand\uh {{{u}_h}}

\newcommand\buh {{\bar{u}_h}}

\newcommand\uhk {{{u}_h^{(k)}}}

\newcommand\uhN {{{u}_h^{(0)}}}
\newcommand\uk {{\tilde{u}^{(k)}}}
\newcommand\uhkM {{{u}_h^{(k-1)}}}
\newcommand\zhk {{{z}_h^{(k)}}}

\newcommand\zk {{\tilde{z}^{(k)}}}

\newcommand\uhl {{{u}_{h,\ell}}}
\newcommand\zhl {{{z}_{h,\ell}}}
\newcommand\uhlk {{{u}_{h,\ell}^{(k)}}}
\newcommand\zhlk {{{z}_{h,\ell}^{(k)}}}

\newcommand\vh {{{v}_h}} 
\newcommand\zh {{{z}_h}}

\newcommand\TTT {\mathcal{T}}
\newcommand\Th {\TTT_h}
\newcommand\Thl {\TTT_{h,\ell}}

\newcommand\Fh{{{\mathcal F}_h}}
\newcommand\FhI{{{\mathcal F}_h^I}}
\newcommand\FhIB{{{\mathcal F}_h}}  

\newcommand\FhB{{{\mathcal F}_h^B}}

\newcommand\bHk{H^{k}(\Th)}
\newcommand\bHj{H^{1}(\Th)}
\newcommand\bHd{H^{2}(\Th)}

\newcommand{\norm}[2]{\left\|#1\right\|_{#2} }

\newcommand\aaa {{a}}
\newcommand\aaaP {{a^{\prime}}}
\newcommand\ah {{a_h}}
\newcommand\ahP {{a_h^{\prime}}}
\newcommand\ahL {{a_h^{\scriptscriptstyle\mathrm{L}}}}

\newcommand\aht {{\tilde{a}_h}}

\newcommand\om{{\Omega}}
\newcommand\J{J}
\newcommand\Jh{J_h}
\newcommand\Jht{\tilde{J}_h}
\newcommand\JhL{J_h^{\scriptscriptstyle\mathrm{L}}}
\newcommand\JP {{J^{\prime}}}

\newcommand\plusminus{^{\scriptscriptstyle (\pm)}}
\newcommand\plus{^{\scriptscriptstyle (+)}}
\newcommand\minus{^{\scriptscriptstyle (-)}}

\renewcommand\bkn{\vec{n}}
\newcommand\bknK {{\bkn_{\!\!\;\scriptscriptstyle K}}}
\newcommand\bkng {{\bkn_{\gamma}}}

\newcommand\V {V}
\newcommand\Vh {V_h}
\newcommand\VhP {V_h^+}

\newcommand\Whm {W_{h,{\ell}}}

\newcommand\Wh {W_h}
\newcommand\WhP {W_h^+}
\newcommand\W {W}

\providecommand{\res}[2]{r_h(#1)(#2) }

\providecommand{\resD}[2]{r_h^*(#1)(#2) }

\providecommand{\ResNLk}[1]{\rho_h(\uhk)(#1) }
\providecommand{\ResNLDk}[1]{\rho_h^*(\uhk, \zhk)(#1) }

\newcommand\Rtri {\mathcal{R}_h^{(3)} }

\newcommand\etaA {\eta^{\scriptscriptstyle\mathrm{A}}}
\newcommand\etaS {\eta^{\scriptscriptstyle\mathrm{S}}}
\newcommand\etaL {\eta^{\scriptscriptstyle\mathrm{L}}}
\newcommand\etaJ {\eta^{\scriptscriptstyle\mathrm{J}}}

\newcommand\etaI {\eta^{\scriptscriptstyle\mathrm{I}}}
\newcommand\etaK {\eta^{\scriptscriptstyle\mathrm{I}}_K}

\newcommand\estA {e_h^{\scriptscriptstyle\mathrm{A}}}
\newcommand\estS {e_h^{\scriptscriptstyle\mathrm{S}}}
\newcommand\estL {e_h^{\scriptscriptstyle\mathrm{L}}}
\newcommand\estJ {e_h^{\scriptscriptstyle\mathrm{J}}}
\newcommand\estR {e_h^{\scriptscriptstyle\mathrm{R}}}

\newcommand\estLAx {e_h^{\scriptscriptstyle\mathrm{A\star}}}
\newcommand\estLJx {e_h^{\scriptscriptstyle\mathrm{J\star}}}

\newcommand\ehl {{\tilde{e}_{h}}}
\newcommand\ehlD {{\tilde{e}_{h}^*}}

\newcommand\TOL {\mathrm{TOL}} 
\newcommand\jump[1]{[\![{#1}]\!]}
\newcommand\aver[1]{\{\!\!\{{#1}\}\!\!\}}
\newcommand{\Lsp}[2]{({#1},{#2})_{\Omega}}

\newcommand\CA{{C_{{A}}}}
\newcommand\CL{{C_{I}}}

\newcommand\DoF{\mathrm{DoF}}

\newcommand\HGhp{{$hp$-HG}}
\newcommand\HGh{{$h$-HG}}
\newcommand\AMAh{{$h$-AMA}}
\newcommand\AMAhp{{$hp$-AMA}}

\newcommand\pK {p_K}

\newcommand\ralg {\bkr}
\newcommand\ieff {i_{\mathrm{eff}}}

\newcommand\Wuhk{P\uhk}
\newcommand\Wzhk{P\zhk}

\usepackage[usenames]{color}

\newcommand{\Smaz}[1]{{}}

\newcommand{\refA}[1]{{{#1}}}
\newcommand{\refB}[1]{{{#1}}}

\begin{document}
\begin{frontmatter}



\title{Goal-oriented error analysis of iterative Galerkin discretizations for nonlinear problems including linearization and algebraic errors\tnoteref{label1}}
\tnotetext[label1]{This work was supported by grant No. 20-01074S of the Czech Science Foundation.}

\author[Prague]{V{\'\i}t Dolej{\v s}{\'\i}}
\ead{dolejsi@karlin.mff.cuni.cz}

\author[Prague]{Scott Congreve}
\ead{congreve@karlin.mff.cuni.cz}

\address[Prague]{Charles University, Faculty of Mathematics and Physics,
Sokolovsk\'a 83, 186 75 Praha, 
Czech Republic}

\begin{abstract}
  We consider the goal-oriented error estimates for
    a linearized iterative solver for nonlinear partial differential
    equations. For the adjoint problem and iterative solver we consider,
    instead of the differentiation of the primal problem, a
    suitable linearization which guarantees the adjoint consistency of the numerical scheme.
    We derive error estimates and
    develop an efficient adaptive algorithm which balances the errors arising
    from the discretization and use of iterative solvers.
    Several numerical examples demonstrate the efficiency of this algorithm.
\end{abstract}
\begin{keyword}
    goal-oriented error estimates \sep nonlinear problems\sep algebraic errors\sep adaptive solvers
    \sep stopping criteria
    
    \MSC 65N30 \sep  65N15 \sep 65N50
\end{keyword}
\end{frontmatter}
  

\section{Introduction}

When computing a numerical approximation to a partial differential equation we are often interested
in the value of a certain solution-dependent target functional rather than the global approximate solution.
This has led, in recent decades, to the development of the goal-oriented error estimates and mesh adaptation
techniques; cf. \cite{RannacherBook,BeckerRannacher01,GileSuli02,FidkowskiDarmofal_AIAA11}
and the references cited therein.
In order to estimate the error of the quantity of interest, an adjoint problem is formulated
and its solution is employed in the error estimates.
While this technique is well-developed for linear problems, for nonlinear problems
a suitable linearization of the primal problem is required; 
see the seminal works \cite{RannacherBook,BeckerRannacher01,GileSuli02} and
some fluid dynamics applications, e.g., in
\cite{HH06:SIPG2,LoseilleDervieuxAlauzet_JCP10,BalanWoopenMay16,FidkowskiDarmofal_AIAA11}.
The linearization and the setting of the adjoint problem should be formulated such that, \refA{asymptotically,
{\em adjoint consistency} is achieved}; cf., \cite{Hartmann2007Adjoint}.

In the above cited works, the primal problem is linearized by differentiation
with respect to the approximate solution, which aligns
with solving the primal problem by the Newton method since the corresponding Jacobian is available.
Nevertheless, \refA{the primal form is not differentiable for some type of problems.}
Therefore, alternative nonlinear solvers have to be employed,
e.g., Newton-like methods in \cite{DGM-book} or
L-scheme in \cite{RaduALL_JCAM15,ListRadu_ComputGeosc16}. For numerical analysis dealing
with the linearization and iterative solvers, we refer to
\cite{ChaillouSuri_CMAME06,ErnVohralik_SISC13,CongreveWihler_JCAM17,HeidWihler_CAL20}.
In \cite{GO_nonlinear}, we proposed a heuristic framework of goal-oriented error estimates
where the adjoint problem was built on a linearization employed in nonlinear algebraic solvers.
In this article, we present a deeper abstract analysis of this approach
and derive error estimates taking into account the errors arising from
the linearization.
The validity of these estimates is supported by several numerical experiments. This is
the main novelty of this work.

Furthermore, we note that an efficient numerical computation requires
a balance between the discretization error and
errors arising from the inaccurate solution of algebraic systems.
The effect of algebraic errors in the goal-oriented error estimates
has been studied  for the first time in \cite{Rannacher201323} and 
further developed, e.g., 
in
\cite{DiStolfoALL_JNM19,EndtmayerLangerWisk_SISC20,MallikVohralikYousef_JCAM20}.
Extending the results from \cite{GO_nonlinear}, we propose an adaptive method
for the solution of the algebraic systems combining the nonlinear and linear solvers together
with different mesh adaptation techniques. Its performance is demonstrated by several
numerical experiments ranging from simple benchmarks to a practically motivated example. 

The outline of this article is as follows. In Section~\ref{sec:DWR}, we derive error estimates
for an iterative solver for nonlinear problems based on a linearization of the primal problem.
We show in Section~\ref{sec:exam}, for a particular example, the formulation and linearization
in both continuous and discontinuous Galerkin finite element methods.
In Section~\ref{sec:adapt}, we discuss an efficient adaptive algorithm. Then, we demonstrate in Section~\ref{sec:numer},
via numerical experiments, the efficiency of the method and accuracy of the error estimate.
Finally, in Section~\ref{sec:concl}, we summarize the results presented in this article.

\section{Goal-oriented error estimates}
\label{sec:DWR}

Given a semilinear form  $\aaa:W\times V \to \R$ associated to the variational formulation
of a nonlinear problem, where $V$ and $W$ are suitable functional spaces, the weak solution $\u\in W$ is given by
\begin{align}
  \label{PP}
  \aaa(\u, \vp) = 0\qquad \forall\vp\in V.
\end{align}
The boundary conditions are realized either by the choice of $V$ and $W$, or they are directly included
in the form $a$.

The approximate solution of \eqref{PP} is sought in the finite dimensional space $\Wh$ 
which consists of piecewise polynomial functions on a mesh $\Th$, with a finite dimensional test space $\Vh$.
We admit $\Vh\subset V$ and $\Wh\subset W$ as well as $\Vh\not\subset V$ and $\Wh\not\subset W$. Moreover, let $\V(h)$ be a functional
space such that $V\subset \V(h)$ and $\Vh\subset \V(h)$, and similarly let $\W(h)$ be a functional space such that $W\subset \W(h)$ and $\Wh\subset \W(h)$. For conforming finite element methods, we
set $\V(h)\coloneqq V$ and $\W(h)\coloneqq W$. For discontinuous Galerkin methods, $\V(h)$ and $\W(h)$ are broken Sobolev spaces,
see \eqref{bss}. 
The quantity of interest is given by a possibly nonlinear functional $\J:\W(h)\to\R$.

Let $\ah:\W(h)\times\V(h)\to\R$ be a semilinear
form representing the discretization of $\aaa$ by a suitable numerical method.
Then, $\uh\in \Wh$ is the approximate solution of \eqref{PP} if
\begin{align}
    \label{PPh}
    \ah(\uh, \vph) = 0\qquad \forall\vph\in \Vh.
\end{align}
The problem \eqref{PPh} represents a system of nonlinear algebraic equations which has to be
solved iteratively. Hence, only an approximation of $\uh$ is available. We assume
that $\ah$ is consistent, i.e., if $\u$ is the solution of \eqref{PP} then
\begin{align}
    \label{cons}
    \ah(\u, \vp) = 0\qquad \forall\vp\in \V(h).
\end{align}

\subsection{Differentiation based goal-oriented error estimates}
\label{sec:DWRa}

We now briefly recall the general framework of the goal-oriented error estimates
according to \cite{Rannacher201323} for an abstract nonlinear problem which is based on the
differentiation of the primal problem.
By $\JP[{\vh}](\vp)$ we denote the Fr\'echet derivative of $J$ at ${\vh}$
along the direction
$\vp$. Similarly, $\aaaP[{v}](\vp, \cdot)$ and $\ahP[{\vh}](\vp, \cdot)$ denote
the derivative of $\aaa$ and $\ah$ with respect to its
first argument at ${v}$ and ${\vh}$ along the direction $\vp$, respectively.
{We assume that $\aaaP[{v}](\cdot, \cdot):\bar{W}\times\bar{V}\to \IR$ and
  $\ahP[{\vh}](\cdot, \cdot):\W(h)\times\V(h)\to\R$
  where the spaces $\bar{V}$ and $\bar{W}$ differ from $V$ and $W$ in general, depending on particular problems.} 
Then, the adjoint problem (linearized at ${v}\in\W(h)$) reads:
find $\z\in \bar{V}$ such that
\begin{align}
  \label{i3D}
  \aaaP[{v}] (\vp, \z) = \JP[{v}](\vp) \qquad \forall \vp \in \bar{W}.
\end{align}
The discrete adjoint problem corresponding to \eqref{i3D} (linearized at $\uh\in\W(h)$)
is formulated as: find $\zh\in \Vh$ such that
\begin{align}
  \label{i3}
  \ahP[\uh] (\vph, \zh) = \JP[\uh](\vph) \qquad \forall \vp \in \Wh.
\end{align}

\begin{theorem}\cite[Proposition~3.1]{Rannacher201323}
  \label{thm:rannacher}
  Let $\uhk\in \Wh$ and $\zhk \in \Vh$ be the inexact solutions of the problems
  \eqref{PPh} and \eqref{i3}, respectively.
  Then the error of the quantity of interest satisfies
  \begin{align} 
    \label{aEE}
    J(\u) - J(\uhk) = & \tfrac{1}{2} \ResNLk{\z - \zhk}  + \tfrac{1}{2} \ResNLDk{\u - \uhk} 
     + \ResNLk{\zhk} + \Rtri, 
  \end{align}
  where $ \ResNLk{\psi}  := -  \ah(\uhk, \psi)$ and 
  $\ResNLDk{\psi} := \JP[\uhk](\psi) -  \ahP[\uhk] (\psi, \zhk), \ \psi\in \W(h)$ denote
  the residuals of the primal and adjoint problems, respectively, and
  setting $\ehl = \u-\uhk$ and $\ehlD = \z-\zhk$,   the reminder term is given by \\
$    \Rtri =\frac12\int\nolimits_0^1\Big\{  J^{\prime\prime\prime}[\uhk+s\ehl](\ehl,\ehl,\ehl)
    - a_h^{\prime\prime\prime}[\uhk+s\ehl](\ehl,\ehl,\ehl,\zhk+s\ehlD)
    - 3 a_h^{\prime\prime}[\uhk+s\ehl](\ehl,\ehl,\ehlD) {\Big \}s(s-1)}
    \ds $.
\end{theorem}

We note that the remainder term $\Rtri$ is of order $O(|\ehl|^3, |\ehlD|^3)$
and it is usually neglected,
an exception is the numerical analysis in \cite{EndtmayerLangerWisk_SISC20}.
The estimate \eqref{aEE} depends on higher order derivatives
$a_h^{\prime\prime\prime}$ and $J^{\prime\prime\prime}$ and, hence, a sufficient regularity of
$\ah$ and $\J$ is required.

\subsection{Adjoint problem based on the linearization of $\ah$}
We now consider, instead, the generalization of the setting of the adjoint problem considered in 
\cite{GO_nonlinear}. Here, rather than employ the differentiation of forms $\ah$ and $\J$ we instead study their linearization
as used in the iterative solvers.
The relation \eqref{PPh} represents the system of nonlinear algebraic equations
which is solved by an iterative method based on a suitable linearization of $\ah$.
Therefore, we assume that there exist forms $\ahL:\W(h)\times\W(h)\times\V(h)\to\IR$ and
$\aht:\W(h)\times\V(h)\to\IR$ which are consistent with $\ah$ by 
\begin{align}
  \label{linP}
  \ah(\uh,\vph) = \ahL( \uh, \uh, \vph) -  \aht(\uh,\vph)\qquad \forall
  \uh\in\W(h)\,\forall\vph\in\V(h),
\end{align}
where $\ahL$ is linear in its second and third arguments, and $\aht$
is linear in its second argument. Particular examples of the forms $\ahL$ and $\aht$
are given in Section~\ref{sec:exam}.

Using \eqref{linP},  we define the iterative process for the solution of \eqref{PPh}.
Let $\uhN\in\Wh$ be an initial approximation of $\uh$, 
we set the sequence $\uhk\in\Wh,\ k=1,2,\dots$ by
\begin{align}
  \label{newtonP}
  \ahL( \uhkM, \uhk, \vph) =  \aht(\uhkM,\vph)\quad \forall\vph\in\Vh,\ k=1,2,\dots.
\end{align}
This identity exhibits a system of linear algebraic equations which have to be solved by a suitable
solver. We note that in order to improve the convergence $\uhk\to\uh$, a damping factor
has to be included, see, e.g., \cite[Section~8.4.4]{DGM-book}.

We define the iterative primal residual as
\begin{align}
  \label{res}
\res{\vh,\uh}{\vph} := - \ahL(\vh,\uh,\vph) + \aht(\vh,\vph),\quad \uh,\vh\in\Wh, \vph\in\Vh.
\end{align}
Obviously, if the first two arguments to $r_h$ are the same function then $r_h$ is equivalent to $\rho_h$ from Theorem~\ref{thm:rannacher}.
\refB{As it is possible to select $\aht$ for almost any choice of $\ahL$ such that \eqref{linP} is satisfied, then \eqref{newtonP} covers most nonlinear iterative techniques; for example, Ka\v{c}anov or Zarantonello. Zarantonello iterations, cf. \cite[Section 25.4]{ZeidlerIIB}, for example is given by $\ahL(\uh,\cdot,\vh) := \delta^{-1}(\uh,\vh)_X$, for some inner product $X$ and constant $\delta$, with $\aht(\uh,\vh)=\delta^{-1}(\uh,\vh)_X-\ah(\uh,\vh)$.
  Section~\ref{sec:exam} gives concrete examples for Ka\v{c}anov iterations.}
Moreover, \eqref{newtonP} covers also the Newton method
provided that we set $\ahL(\uh,\cdot,\cdot) := \ahP[\uh] (\cdot,\cdot)$.
Therefore, this setting is more general
than the one from Section~\ref{sec:DWRa}.

In has been shown in \cite{Hartmann2007Adjoint} that in order to guarantee 
the adjoint {consistency} of the method,  modification of the target functional
$\J$ is required even for linear problems.
Therefore, we introduce a new functional $\Jh:\W(h)\to\R$ which is consistent with
$\J$ such that $\Jh(\u) = \J(\u)$,
where $u$ is the solution of \eqref{PP}.
Moreover, we assume that there exists a linearization of $\Jh$,
namely the forms $\JhL:\W(h)\times\W(h)\to\R$ and  $\Jht:\W(h)\to\R$ which are
consistent with $\Jh$ by 
\begin{align} 
\label{JhL}
\Jh(\vph) = \JhL(\vph, \vph) + \Jht(\vph)\qquad \forall \vph\in\Wh,
\end{align}
where form $\JhL$ is linear in its second argument. Form $\Jht$ is often independent of $\vph$,
e.g., if it arises from the replacement of $\J$ by $\Jh$.
In this case form $\Jht$ does not influence the error since we have
$\Jh(\u) - \Jh(\uh) = \JhL(\u, \u) - \JhL(\uh, \uh)$.

We introduce the adjoint problem using the linearized forms $\ahL$ and $\JhL$.
We say that $\zh \in \Vh$ is the discrete adjoint solution of the adjoint problem
(linearized at $\uh\in\Wh$) if it satisfies
\begin{align} 
  \label{DPh}
  \ahL( \uh, \vph, \zh) = \JhL(\uh,\vph) \qquad \forall \vph \in \Wh,
\end{align}
where $\ahL$ and $\JhL$ are the forms from \eqref{linP} and \eqref{JhL}, respectively.
The corresponding adjoint residual is given by
\begin{align}
  \label{resD} 
  \resD{\uh,\zh}{\vph} :=  \JhL(\uh,\vph) - \ahL( \uh, \vph, \zh),\quad \uh,\vph\in\W(h),\
  \zh\in\V(h).
\end{align}
At each step of the iterative process \eqref{newtonP},
we can define the corresponding adjoint approximation by finding $\zhk\in\Vh$ such that
\begin{align}
    \label{iteradjoint}
    \ahL( \uhkM, \psih,\zhk) =  \JhL(\uhkM,\psih)\qquad \forall\psih\in\Wh.
\end{align}

Finally, we assume that discretization \eqref{PPh} is adjoint consistent;
i.e, if $u$ is the weak solution of \eqref{PP} and $z\in V$
is the weak solution of the adjoint problem then
\begin{align}
  \label{consD}
  \ahL( \u, \vph, \z) = \JhL(\u, \vph)  \qquad \forall\vph\in\W (h).
\end{align}

\subsection{Error estimates based on the linearization of $\ah$ and $\Jh$}\label{sec:estim2}
We now derive the goal-oriented error estimation of
the numerical solution obtained at each step $k=0,1,\dots$
of the iterative process given by \eqref{newtonP}.
We introduce the auxiliary primal problem of finding $\uk\in\W(h)$ such that
\begin{align}
    \label{recon_u}
    \ahL( \uhkM, \uk, \vp) =  \aht(\uhkM,\vp)\qquad \forall\vp\in\V(h),
\end{align}
and the auxiliary adjoint problem of finding $\zk\in\V(h)$ such that
\begin{align}
    \label{recon_z}
    \ahL( \uhkM, \psi,\zk) =  \JhL(\uhkM,\psi)\qquad \forall\psi\in\W(h).
\end{align}
We note that $\uk\in\W(h)$ and $\zk\in\V(h)$ are \emph{reconstructions}
(cf. \cite{MakridakisNochetto2003}) in the sense that $\uhk\in\Wh$ from \eqref{newtonP}
and $\zhk\in\Vh$ from \eqref{iteradjoint} are the Galerkin approximations of $\uk$ and $\zk$,
respectively.

\begin{theorem}
  \label{thm:main}
  Let $\uhk\in\Wh$ and $\zhk\in\Vh$ be the numerical approximations given by \eqref{newtonP} and \eqref{iteradjoint}, respectively. Then, the error of the quantity of interest satisfies
  \begin{align}
    \label{EST1}
    \J(u)-\Jh(\uhk) = \estS + \estA + \estL + \estJ + \estR,
    \end{align}
    where
    \begin{align}
    \label{EST2}
    \estS &= \frac12 \left( \res{\uhkM,\uhk}{\zk-\zhk} + \resD{\uhkM, \zhk}{\uk-\uhk} \right),
    \qquad   \estA = \res{\uhkM,\uhk}{\zhk}, \\
    \estL &= \ahL(\u,\u-\uhk,\z)-\ahL(\uhkM,\uk-\uhk,\zk),
    \qquad \estJ=  \JhL(\u,\uhk)-\JhL(\uhk,\uhk),  
    \qquad 
    \estR = \Jht(\u)-\Jht(\uhk). \notag
    \end{align}
    Moreover, the linearization part of the error can be written by an alternative formula
    \begin{align}
      \label{EST3}
      \estL +\estJ &=\estLAx + \estLJx \\
      &\coloneqq \left(\ahL(\u,\u,\z-\zhk)-\ahL(\uhkM,\uk,\zk-\zhk)\right)+ 
        \left(\aht(\u,\zhk)-\aht(\uhkM,\zhk)+\JhL(\uhkM,\uhk) - \JhL(\uhk,\uhk)\right). \notag
    \end{align}
\end{theorem}
\begin{remark}
  \label{rem:null}
  If $\Jh$ is linear, then $\JhL$ is independent of its first argument; therefore, with a linear quantity of interest the $\estJ$ term disappears.
  Similarly, if $\aht$ is independent of its first argument, then
  $\estLJx=\JhL(\uhkM,\uhk) - \JhL(\uhk,\uhk)$.
  Furthermore, $\Jht$ is often independent of its arguments; hence, $\estR=0$.
\end{remark}
\begin{proof}
  From \eqref{JhL}, \eqref{recon_u}, and \eqref{res} we can write the error of the quantity of interest as
  \begin{align}
    \label{error}
    J(\u) - \Jh(\uhk) &= \left(\JhL(\uhkM,\uk) - \JhL(\uhkM,\uhk) - \ahL(\uhkM,\uk-\uhk, \zhk)\right) + \res{\uhkM,\uhk}{\zhk}\\
    &\quad 
    + \left( \JhL(\u,\u) - \JhL(\uhkM,\uk) + \JhL(\uhkM,\uhk) - \JhL(\uhk,\uhk) \right) 
    + \left( \Jht(\u)-\Jht(\uhk)\right) \notag\\
    &\coloneqq\xi_1 + \estA + \xi_2 + \estR. \notag
\end{align}
For $\xi_1$, by \eqref{recon_z}, we have 
\begin{align}
    \label{xi1_initial}
    \xi_1 & =  \JhL(\uhkM,\uk) - \JhL(\uhkM,\uhk)  - \ahL(\uhkM,\uk-\uhk, \zhk) \\
    &= \ahL(\uhkM, \uk,\zk) - \ahL(\uhkM, \uhk,\zk) - \ahL(\uhkM,\uk-\uhk, \zhk). \notag
\end{align}
By applying \eqref{recon_u} we have 
\begin{align}
\label{xi1_1}
\xi_1 &=  \aht(\uhkM,\zk) - \ahL(\uhkM, \uhk,\zk) - \aht(\uhkM,\zhk) + \ahL(\uhkM, \uhk,\zhk)  \\
& =  -\ahL(\uhkM, \uhk,\zk-\zhk) +\aht(\uhkM,\zk-\zhk) \eqqcolon \res{\uhkM,\uhk}{\zk-\zhk}. \notag
\end{align}
Furthermore, from \eqref{xi1_initial} we have
\begin{align}
  \label{xi1_2}
  \xi_1 &=  \ahL(\uhkM, \uk-\uhk,\zk)-\ahL(\uhkM, \uk-\uhk,\zhk)\\
  & =  \JhL(\uhkM, \uk-\uhk)-\ahL(\uhkM, \uk-\uhk,\zhk)  \eqqcolon \resD{\uhkM,\zhk}{\uk-\uhk}.
  \notag 
\end{align}
Combining \eqref{xi1_1} and \eqref{xi1_2} we get that
\begin{align}
    \label{xi1}
    \xi_1 = \tfrac12 \left( \res{\uhkM,\uhk}{\zk-\zhk} + \resD{\uhkM, \zhk}{\uk-\uhk} \right) \eqqcolon \estS.
\end{align}

Moreover, from \eqref{consD} -- \eqref{recon_z}, we have
\begin{align}
  \label{X4}
  & \JhL(\u,\u) - \JhL(\uhkM,\uk)  = \ahL(\u,\u,\z) - \ahL(\uhkM, \uk,\zk) \\
  &\qquad= \ahL(\u,\u,\z) - \ahL(\uhkM, \uk,\zk) - ( \ahL(\u,\uhk,\z) - \JhL(\u,\uhk))
  + ( \ahL(\uhkM,\uhk,\zk) - \JhL(\uhkM,\uhk)) \notag \\
  &\qquad= \left(\ahL(\u,\u-\uhk,\z) - \ahL(\uhkM, \uk-\uhk,\zk)\right) +
  \left( \JhL(\u,\uhk) - \JhL(\uhkM,\uhk)\right), \notag
\end{align}
which together with the definition of $\xi_2$ in \eqref{error} gives $\xi_2= \estA + \estJ$;
cf. \eqref{EST2}.
Similarly, we derive the alternative formula \eqref{EST3}.
Again, using \eqref{consD} -- \eqref{recon_z},
we obtain
\begin{align}
  \label{X5}
  & \JhL(\u,\u) - \JhL(\uhkM,\uk)  = \ahL(\u,\u,\z) - \ahL(\uhkM, \uk,\zk) \\
  &\qquad= \ahL(\u,\u,\z) - \ahL(\uhkM, \uk,\zk) - ( \ahL(\u,\u,\zhk) - \aht(\u,\zhk))
  + ( \ahL(\uhkM,\uk,\zhk) - \aht(\uhkM,\zhk))  \notag \\
  &\qquad= \left( \ahL(\u,\u,\z-\zhk) - \ahL(\uhkM, \uk,\zk-\zhk) \right)
  + \left( \aht(\u,\zhk) - \aht(\uhkM,\zhk) \right). \notag
\end{align}
Hence, the definition of $\xi_2$ in \eqref{error} and \eqref{X5} gives
$  \xi_2 =  \estLAx + \estLJx$. \qedhere
\end{proof}

\subsection{\refB{Computable error estimates}}\label{sec:comp_err}

\refB{The main theoretical error estimate \eqref{EST1}--\eqref{EST2} formulated in Theorem~\ref{thm:main}
contains the exact solutions $\u$ and $\z$ of the primal and adjoint problems, respectively,
as well as the linearized solutions $\uk$ and $\zk$, which 
are not available and they have to be approximated.
One possibility, is to construct 
higher order approximations from the available approximate solutions $\uhk$ and $\zhk$
of the discrete problems.
While the approximate solutions $\uhk$ and $\zhk$ are sought in the space $\Wh$, the reconstructions
must belong to a rich space denoted $\WhP$. We then define a reconstruction operator $\krr:\Wh\to\WhP$.}

\refB{Therefore, we approximate the unknown functions in \eqref{EST2} as
\begin{align}
	\label{CE5}
	\u\approx \Krr{\uhk},\qquad
	\uk\approx \Krr{\uhk},\qquad
	\z\approx \Krr{\zhk},\qquad
	\zk\approx \Krr{\zhk}.
\end{align}
Both $\u$ and $\uk$ are approximated by the same function since $\uhk$
is the only available information. 
The presented numerical experiments in Section~\ref{sec:numer} show that these approximations
give a reasonable computational performance.
Finally, in virtue of \eqref{EST1}--\eqref{EST2} and \eqref{CE5},
we define a computable approximation of the error by
\begin{align}
	\label{CE10}
	\Jh(\uhk) - \J(\u) & \approx  \etaI (\uhk,\zhk)
	:= \etaS (\uhk,\zhk) +\etaA (\uhk,\zhk)+\etaJ (\uhk,\zhk)+\etaL (\uhk,\zhk),
\end{align}
where, for simplicity, we omit the explicit dependence on $\uhkM$ and
\begin{align}
	\label{CE10a}
	\etaS (\uhk,\zhk) &:= \frac{1}{2}\Big(\res{\uhkM, \uhk}{\Krr{\zhk}-\zhk}
	+ \resD{\uhk,\zhk}{\Krr{\uhk}-\uhk} \Big),  
	\qquad \etaA (\uhk,\zhk) := \res{\uhkM, \uhk}{\zhk},  \notag \\
	\etaL (\uhk,\zhk) &:=  \ahL\Big(\Krr{\uhk}, \Krr{\uhk}-\uhk,\zhk\Big)
	-  \ahL\Big(  \uhkM, \Krr{\uhk}-\uhk,\zhk\Big), \quad
	\etaJ:= \JhL\Big(\Krr{\uhk},\uhk\Big)-\JhL\Big(\uhk,\uhk\Big). 
\end{align}
Estimator $\etaS$ corresponds to the weighted residual error, $\etaA$ to the algebraic error
and finally, $\etaL$ and $\etaJ$ are the error estimators
arising from the linearization of $\ah$ and $\Jh$,
respectively.
In virtue of Remark~\ref{rem:null}, the linearization estimator $\etaL$ can
be replaced by the alternative formula following from relation \eqref{EST3}.
Moreover, we do not consider the term $\estR$ since it vanishes in our examples.}

\refB{For the purpose of the mesh adaptation, the error estimate \eqref{CE10} has to be localized;
i.e., define local estimates $\etaK(\uhk,\zhk)$, $K\in\Th$, such that
$\etaI (\uhk,\zhk) = \sum_{K\in\Th} \etaK (\uhk,\zhk)$.
The construction is done usually by a partition of unity, cf. \cite{RichterWick_JCAM15};
for details we refer
to \cite[Chapter~7]{AMA-book}.}

\section{Several particular examples}
\label{sec:exam}

In this section, we present the discretization of a concrete problem,
by both the {\em continuous Galerkin} and the {\em symmetric interior penalty Galerkin} (SIPG) variant of the 
{\em discontinuous Galerkin} (DGM) finite element methods. Furthermore, we demonstrate the linearization for the adjoint problem \eqref{DPh} and the iterative process \eqref{newtonP}.

Let $\Om\subset \R^d,\ d=2,3$ be a bounded domain with Lipschitz boundary $\Gamma:=\pd\Om$. Then, we consider the nonlinear diffusion-reaction problem: find $u:\Om\to\R$ such that
\begin{alignat}{2}
    \label{diff_reac}
     -\nabla\cdot \left(\mu(|\nabla u|)\nabla u \right) + d(u)u &= f &&
    \quad\mbox{in }\Omega, \\
    u &= g && \quad\mbox{on }\gom,\notag
\end{alignat}
where $f\in L^2(\Omega)$, $g$ is the trace of a function in $H^1(\Omega)$,
$\mu$ is a strongly monotone and Lipschitz continuous nonlinear function, and $d(u) \in C(\R)$.

Defining the space $W=H^1_g(\Om)\coloneqq\{\vp\in H^1(\Om) : \vp=g \text{ on }\gom \}$ and the space $V=H^1_0(\Om)$ as the space of functions in $H^1(\Omega)$ with zero trace on $\gom$; then, the weak solution $u\in H^1_g(\Om)$ to \eqref{diff_reac} is given by \eqref{PP}, where
\begin{align}
\label{ex_a}
a(u,\varphi) \coloneqq \int_\Om \mu(|\nabla u|)\nabla u \cdot \nabla\vp \dx +  \int_\Om d(u)u\vp \dx - \int_\Om f\vp \dx, \qquad u\in H^1_g(\Om),\vp\in H^1_0(\Om).
\end{align}

In Section~\ref{sec:numer}, we consider the target functional representing the energy associated to the diffusion part of
\eqref{diff_reac} given by $J(u) = \int_\Om \chi_M\, \mu(|\nabla u|) |\nabla u |^2\dx$ where
$\chi_M$ is the characteristic function of a subdomain $M\subset\Om$.
Then the linearization \eqref{JhL} is defined by
\begin{align} 
\label{JhLA}
\JhL(\uh, \vph)  \coloneqq  \int_\Om \chi_M\, \mu(|\nabla \uh|) \nabla \uh \cdot\nabla\vph \dx,
\qquad  \Jht(\vph)= 0.
\end{align}

\subsection{Continuous Galerkin method}
We first consider the formulation and linearization using a continuous Galerkin finite element method. To this end, we let $\Th$ be a regular and shape-regular mesh that partitions $\Om$ into open disjoint simplices $K$ such that $\oO=\bigcup_{K\in\Th} \overline{K}$. For a fixed polynomial degree $p\geq 1$ we introduce the finite element spaces
\begin{align}
    \SShp \coloneqq \{ v_h\in H^1(\Om);\  v_h\vert_K \in P^p(K) \  \forall K\in\Th \} \subset H^1(\Omega), \label{cg_Vh}
\end{align}
$\Vh=\SShp\cap H^1_0(\Omega)$, and $\Wh=\SShp\cap H^1_g(\Omega)$, where $P^p(K)$ is the space of polynomials of total degree at most $p$ on $K$.
As $\Wh\subset H^1_g(\Om)=W$ and $\Vh\subset H^1_0(\Om)=V$, we can set $\W(h)=H^1_g(\Om)$, $\V(h)=H^1_0(\Om)$, and define $a_h$ from \eqref{PPh} as \eqref{linP} 
where
\begin{align}
    \label{cg_ahL}
    \ahL(\buh, \uh, \vph) &\coloneqq \int_\Om \mu(|\nabla \buh|)\nabla \uh \cdot \nabla\vph \dx + \int_\Om d(\buh)\uh\vph \dx,&& \buh,\uh\in H^1_g(\Om),\vph\in H^1_0(\Om), \\
    \label{cg_aht}
    \aht(\cdot,\vph)  &\coloneqq\int_\Om f\vph \dx, && \vph\in H^1_0(\Om).
\end{align}
The approximate solution $\uh\in\Wh$ of \eqref{diff_reac} is, therefore, given by \eqref{PPh}, and
the inexact iterative solution is given, for an initial approximation $\uhN\in\Wh$, by \eqref{newtonP}.
Furthermore, the discrete adjoint solution and its iterative approximation is given by \eqref{DPh}
and \eqref{iteradjoint}, respectively.
\refB{\begin{remark}
		We note that here we have assumed that the Dirichlet boundary condition $g$ belongs to the space $S_h^p$; for example, a constant boundary condition. For more complicated boundary conditions, $g$ must be approximated in $S_h^p$ leading to additional, potentially lower order, error terms.
\end{remark}}

\subsection{Symmetric interior penalty discontinuous Galerkin method}
\label{sec:DGM}
We now consider the formulation and linearization for a discontinuous Galerkin finite element method.
Here, we allow the mesh $\Th$ to contain hanging nodes.
We denote by $\pd K$ and $\bknK$ the boundary of element $K\in\Th$ and the unit outer
normal to $\pd K$, respectively. We also introduce the
{\em broken Sobolev spaces} 
\begin{align}
    \label{bss}
    \bHk := \{\vp\in L^2(\Om);\  \vp|_K\in H^k(K)\ \forall K\in\Th\},\quad k\in\matN.
\end{align}

We denote by $\FhI$ and $\FhB$ the set of all interior faces/edges and boundary faces/edges,
respectively, of the mesh $\Th$. Additionally, we let $\Fh=\FhI\cup\FhB$ denote the set of all
faces/edges in the mesh $\Th$. For each $\gamma\in\Fh$,
we associate the unit normal vector $\bkng$ whose orientation is arbitrary but fixed, and assume $\bkng$ is the outer unit normal to $\Gamma$ for $\gamma\in\FhB$.

Given two adjacent elements, $K\plus$ and $K\minus$, which share an edge $\gamma$, orientated such that $\bkng$ is the outer unit normal with respect to $K\plus$, then we write $\vp|_\gamma\plusminus$ to denote the traces of $\vp\in\bHj$ on $\gamma\in\Fh$, taken from the interior of $K\plusminus$, respectively. We define the mean value and jump on $\gamma\in\FhI$ by
   $\aver{\bkq}_\gamma = \tfrac12\left(\bkq|_{\gamma}\plus + \bkq|_{\gamma}\minus\right)$ and
$\jump{\vp}_\gamma = \left(\vp|_{\gamma}\plus - \vp|_{\gamma}\minus\right)\bkng$, respectively,
for vector functions $\bkq \in[\bHj]^d$ and scalar functions $\vp\in\bHj$.
On a boundary face $\gamma\in\FhB$ we set $\aver{\bkq}_\gamma = \bkq|_{\gamma}\plus$ and $\jump{\vp}_\gamma = \vp|_{\gamma}\plus \bkng$.
For simplicity, we omit the subscript $\gamma$ in $\aver{\cdot}_\gamma$ and $\jump{\cdot}_\gamma$.

The diffusive terms in \eqref{diff_reac} are discretized by the SIPG variant of
DG method according to \cite{DGM-book},
which differs from the technique in \cite{houston-robson-suli}. 
Let $\V(h)=\W(h)\coloneqq\bHd$, then the form $\ah$ from \eqref{PPh} representing the DG discretization of problem \eqref{diff_reac} is given by \eqref{linP} 
where 
\begin{align}
  \label{dg_ahL}
  \ahL(\buh, \uh, \vph)
    &\coloneqq \sum_{K\in\Th} \int_K \mu(|\nabla \buh|)\nabla \uh\cdot\nabla \vph\dx + \sum_{K\in\Th} \int_K d(\buh)\uh\vph\dx
         -\sum_{\gamma\in\FhIB}\int_\gamma\aver{\mu(|\nabla \buh|)\nabla \uh}\cdot\jump{\vph}\ds\\
    &\quad 
        -\sum_{\gamma\in\FhIB}\int_\gamma\aver{\mu(|\nabla \buh|)\nabla \vph}\cdot\jump{\uh}\ds
        +\sum_{\gamma\in\FhIB}\int_\gamma\sigma\jump{\uh}\cdot\jump{\vph}\ds, \notag\\
  \aht(\buh,\vph)
    &\coloneqq \sum_{K\in\Th} \int_K f\vph\dx
     - \sum_{\gamma\in\FhB}\int_\gamma \mu(|\nabla \buh|)\nabla \vph \cdot \bkng g\ds
     + \sum_{\gamma\in\FhB}\int_\gamma \sigma g\vph\ds, \notag
\end{align}
for $\buh,\uh,\vph\in \bHd$,
and $\sigma>0$ is the penalty parameter proportional to the inverse of the diameter of
$\gamma\in\Fh$. We define the discontinuous finite element space
\begin{align}
  \label{dg_Vh}
  \Wh:=\{v_h\in L^2(\Om);\ v_h\vert_K \in P^{p_K}(K),\ K\in \Th\},
\end{align}
where $P^{p_K}(K)$ denotes the space of polynomial functions of total degree at most $p_K$ on $K\in\Th$
and $p_K$ is the local polynomial approximation degree for each $K\in\Th$. We also set $\Vh=\Wh$.

The approximate solution $\uh\in\Wh$ of \eqref{diff_reac} is, therefore, given by \eqref{PPh}, and
the inexact iterative solution is given, for an initial approximation $\uhN\in\Wh$, by \eqref{newtonP}.
Furthermore, the discrete adjoint solution and its iterative approximation is given by \eqref{DPh}
and \eqref{iteradjoint}, respectively. The consistency and adjoint consistency of this
linearization is  derived in \cite{GO_nonlinear}.

\refB{
Given the enriched space
\begin{align}
	\WhP=\VhP:=\{v_h\in L^2(\Om) : v_h\vert_K \in P^{p_K+1}(K),\ K\in \Th\},\label{dg_VhP}
\end{align}
a reconstruction operator $\krr:\Wh\to\WhP$, the
notation $\Wuhk:={\Krr{\uhk}-\uhk}$ and $\Wzhk:={\Krr{\zhk}-\zhk}$,
we can define the computable error bounds \eqref{CE10a} for this formulation with (linearized) target functional \eqref{JhLA} as
\begin{align*}
  & \etaS (\uhk,\zhk) := \frac{1}{2}\Bigg(-\ahL\left(\uhkM,\uhk,\Wzhk\right)
  +\aht\left(\uhkM,\Wzhk\right) 
  +  \JhL\left(\uhk,\Wuhk\right) -\ahL\left(\uhk,\Wuhk,\zhk\right) \Bigg),  \\
  &\ \ = \frac{1}{2}\Bigg(\sum_{K\in\Th} \int_K  \left(f\,\Wzhk 
  - 
  \mu(|\nabla \uhkM|)\nabla \uhk\cdot\nabla\Wzhk
  -
  d(\uhkM)\uhk\,\Wzhk\right) \dx
  \\
  &\qquad
  +
  \sum_{\gamma\in\FhIB}\int_\gamma\aver{\mu(|\nabla \uhkM|)\nabla \uhk}\cdot\jump{\Wzhk}\ds
  +\sum_{\gamma\in\FhI}\int_\gamma
  \left(\aver{\mu(|\nabla \uhkM|)\nabla\Wzhk}-\sigma\jump{\zhk}\right)\cdot\jump{\uhk}\ds
  \\
  &\qquad 
  +\sum_{\gamma\in\FhB}\int_\gamma
  \left(\mu(|\nabla \uhkM|)\nabla\Wzhk\cdot\bkng-\sigma\Wzhk\right)(\uhk-g)\ds
  \Bigg) \\
  &\quad+\frac{1}{2}\Bigg(\sum_{K\in\Th}
  \int_{K} \left( \chi_M\, \mu(|\nabla \uhk|) \nabla \uhk \cdot\nabla\,\Wuhk 
  -
  \mu(|\nabla \uhk|)\nabla\Wuhk\cdot\nabla \zhk 
  - d(\uhk)\,\Wuhk\zhk\right) \dx
  \\
  &\qquad
  +\sum_{\gamma\in\FhIB}\int_\gamma\aver{\mu(|\nabla \uhk|)\nabla\Wuhk}\cdot\jump{\zhk}\ds
  +\sum_{\gamma\in\FhIB}\int_\gamma\left(\aver{\mu(|\nabla \uhk|)\nabla \zhk}-\sigma\jump{\zhk}\right)\cdot\jump{\Wuhk}\ds
  \Bigg), \\
  \\
  & \etaA (\uhk,\zhk) := -\ahL\left(\uhkM, \uhk, \zhk\right) +\aht\left(\uhkM, \zhk\right),  \\
  &= \sum_{K\in\Th} \int_K \left( f\zhk 
  -
  \mu(|\nabla \uhkM|)\nabla \uhk\cdot\nabla \zhk
  - 
  d(\uhkM)\uhk\zhk \right) \dx
  +\sum_{\gamma\in\FhIB}\int_\gamma\aver{\mu(|\nabla \uhkM|)\nabla \uhk}\cdot\jump{\zhk}\ds
  \\
  &\qquad
  +
  \sum_{\gamma\in\FhI}\int_\gamma\left(\aver{\mu(|\nabla \uhkM|)\nabla \zhk} - \sigma\jump{\zhk}
  \right)\cdot\jump{\uhk}\ds
  +\sum_{\gamma\in\FhB}\int_\gamma\left(
  \mu(|\nabla \uhkM|)\nabla \zhk\cdot\bkng -\sigma \zhk \right)(\uhk-g)\ds,
  \\
  \\
  &\etaL (\uhk,\zhk) :=  \ahL\Big(\Krr{\uhk}, \Wuhk,\zhk\Big)
	-  \ahL\Big(  \uhkM, \Wuhk,\zhk\Big) \\    
        &\ \ = \sum_{K\in\Th} \int_K \left(
        (\mu(|\nabla \Krr{\uhk}|)-\mu(|\nabla \uhkM|))\nabla\Wuhk\cdot\nabla \zhk
        + (d(\Krr{\uhk})-d(\uhkM))\,\Wuhk\zhk\right) \dx \\
        &\quad
        -\sum_{\gamma\in\FhIB}\int_\gamma
        \left(
        \aver{(\mu(|\nabla \Krr{\uhk}|)-\mu(|\nabla \uhkM|))\nabla\Wuhk}\cdot\jump{\zhk}
        - \aver{(\mu(|\nabla \Krr{\uhk}|)-\mu(|\nabla \uhkM|))\nabla \zhk}\cdot\jump{\Wuhk}
        \right)\ds, \\
&	\etaJ(\uhk,\zhk) := \JhL\Big(\Krr{\uhk},\uhk\Big)-\JhL\Big(\uhk,\uhk\Big) =  \sum_{K\in\Th} \int_{K} \chi_M\left( \mu(|\nabla \Krr{\uhk}|) \nabla \Krr{\uhk} - \mu(|\nabla \uhk|) \nabla \uhk\right) \cdot\nabla \uhk\dx.
\end{align*}}
\section{Adaptive algorithm}
\label{sec:adapt}

\subsection{\refB{Higher-order reconstruction}} 

\refB{In order to derive computable error bounds for the discontinuous Galerkin numerical experiments performed in Section~\ref{sec:numer} we need to construct a higher-order reconstruction operator $\krr:\Wh\to\WhP$, see \eqref{dg_Vh} and \eqref{dg_VhP}, as mentioned in Section~\ref{sec:comp_err}; cf.~\eqref{CE5}.}
Based on our extensive experience, we employ the least-square reconstruction
technique from \cite[Section~7.1]{ESCO-16}, which
is sufficiently robust even for problems having singularities.
Let $\vh\in\Wh$, for each $K\in\Th$ we define a patch
$D_K$ consisting of triangles sharing at least a vertex with $K$. 
Then, we seek a function $v_K\in P^{\pK+1}(D_K)$ which minimizes $\norm{\vh - v_K}{D_K}$
and set $\Krr{\vh}|_K := v_K|_K,\ K\in\Th$.
For a numerical study of the accuracy, we refer to \cite{hp_reconstr}.
However, we note that the higher-order reconstruction for anisotropic $hp$-meshes appears to
be a weak point of our technique and, hence, requires further research.

\subsection{Adaptive solver for nonlinear algebraic system}

In this section, we shortly describe the solution strategy
of the nonlinear primal problem \eqref{PPh} and the linear adjoint problem \eqref{DPh}.
As mentioned above, we employ the linearization
\eqref{linP} and define a sequence $\uhk\in\Shp,\ k=1,2,\dots$ iteratively by \eqref{newtonP}.
Similarly, the approximate solution of adjoint problem $\zhk\in\Shp,\ k=1,2,\dots$
is defined by \eqref{iteradjoint}.

For each $k=1,2,\dots$, equations \eqref{newtonP} and \eqref{iteradjoint}
represent linear algebraic systems.
Since it is sufficient to solve them approximately,
the use of an iterative solver is preferable.
In \cite{DolTich_BiCG}, we introduced the technique based on the BiCG solver
which admits to solve both systems simultaneously. However, in some situations, the GMRES solver
seems to be more efficient. The iterative solvers are accelerated by the
ILU(0)-block preconditioner.

A very important question is the choice of the stopping criteria for nonlinear as well as
linear solvers, \refB{cf. \cite{ErnVohralik_SISC13,HaberlAll_NM21,HeidAll_CMAM21}.}
In virtue of \eqref{CE10}, we terminate the nonlinear solver
for an iteration $k\ge 1$ such that
\begin{align}
  \label{CE15}
  |\etaA (\uhk,\zhk)| \leq \CA  |\etaS (\uhk,\zhk)|,
\end{align}
where $\CA\in(0,1)$. In practical examples, we use $\CA=0.01$.
However, the evaluation
of $\etaS$ is much more expensive in comparison to $\etaA$ since $\etaS$ requires the
evaluation of the reconstructions $\Krr{\uhk}$ and $\Krr{\zhk}$. An acceleration can be achieved
by evaluating $\etaS$ only for selected iterations $k$. Typically, when condition \eqref{CE15}
is not valid and the next iteration is necessary, we can avoid updating $\etaS$ and the value
from the previous $k$ can be employed.

Concerning the linear iterative solver, \refB{the usual approach is to balance the error arising
  from the linear solver against the linearization and discretization error,
  cf. \cite{ErnVohralik_SISC13, HaberlAll_NM21}.
  However, the evaluation of those criteria requires also some computation time;
  therefore, based on our experience, we use a simpler criterion.
  The idea is to perform only a few iterations of the linear solver
  and then test criterion \eqref{CE15}. Particularly, }
let $\ralg_{k,l}$  and $\ralg_{k,l}^*$, $l=0,1,\dots$  denote the preconditioned residual vectors
of the linear algebraic systems  \eqref{PPh} and \eqref{DPh},
respectively, achieved
in the $l$-th iteration of the linear solver. Then the linear solver is stopped
when
\begin{align}
  \label{CE16}
  |\ralg_{k,l}| \leq \CL |\ralg_{k,0}|\qquad\mbox{and}\qquad
  |\ralg_{k,l}^*| \leq \CL |\ralg_{k,0}^*|,
\end{align}
where $|\cdot|$ is the Euclidean norm and
$\CL\in(0,1)$ is a suitable constant, typical value is $\CL=0.01$.
This means that the solver is stopped when the preconditioned residual is decreased
by factor 100 (in comparison to the initial residual). This condition may seem to be weak but
we need only an approximate solution of the linear algebraic systems.
\refB{We note that vectors $\ralg_{k,l}$, $l=0,1,\dots$ are automatically available in
 the iterative solvers so no additional computation time is required.}

\subsection{Adaptive mesh algorithm}

The goal of the adaptive mesh algorithm is to obtain a finite element space $\Wh$ and
the corresponding value of the quantity of interest $J(\uh)$, for $\uh\in\Wh$ being the approximate
solution of the primal problem, such that
\begin{align}
  \label{CD1}
  |\J(\u) - \J(\uh)|\approx |\etaI (\uh,\zh) | \leq \TOL,
\end{align}
where $\zh$ is the approximate solution of the adjoint problem,
$\etaI$ is the error estimate \eqref{CE10} and $\TOL>0$ is the given tolerance.
This problem is solved iteratively by defining a sequence of
meshes $\Thl$, spaces $\Shpl$, $\ell=0,1,\dots$
and the corresponding approximations $\uhl,\zhl\in\Shpl$.
If condition \eqref{CD1} is not satisfied,
we adapt the finite element space (and the corresponding mesh)
by one of the techniques described below.

\subsubsection{Anisotropic $hp$-mesh adaptation method}
\label{sec:AMA}
This approach is based on a complete re-meshing of the computational domain and 
the corresponding polynomial approximation degrees. For the detailed method description
we refer to the monograph \cite[Section~7.4]{AMA-book}; here, we mention only the main idea.
In the same manner as in \cite{RannacherBook,BeckerRannacher01}, we apply the discrete
and continuous Cauchy inequalities and re-write the estimate \eqref{CE10}--\eqref{CE10a} as
the sum of several residuals multiplied by weights (= interpolation errors of the
reconstructed primal or adjoint solutions).
Then, we optimize the shape of elements and the polynomial approximation degrees
in such a way that we minimize these weights and the residuals are kept fixed.
We only mention that the theoretical as well as practical results are based on
the so-called continuous  mesh and error models;
cf. \cite{LoseilleAlauzet11a,ESCO-16}. Hereafter, we denote
this method as {\AMAhp}. In the case that we 
keep polynomial degree fixed for all mesh elements, we have the $h$-variant of this method denoted as {\AMAh}.

\subsubsection{Isotropic refinement with hanging nodes}
\label{sec:HG}

For a comparison, we employ also a standard refinement method, where at each adaptation
level, we mark a fixed ratio of elements having the largest value of $|\etaK|$. Typically,
we mark 10\% of elements. Then each marked
triangle is split onto 4 similar sub-elements, i.e., hanging nodes arise. This method is denoted
as {\HGh}.

Moreover, we use the $hp$-variant of this technique, where the
polynomial degree of marked elements can be increased instead of $h$-refinement.
We use a similar criterion as the {\AMAhp} method (cf.~ Section~\ref{sec:AMA}) 
and denote this technique as {\HGhp}.
Note, this method has not been fully developed and will be the subject
of further research.

\section{Numerical experiments}
\label{sec:numer}

We present several numerical examples demonstrating the accuracy of the error estimator
and the performance of the adaptive algorithm described in Section~\ref{sec:adapt}
with discontinuous Galerkin method, cf.~\ref{sec:DGM}.
Some of the examples are benchmarks appearing in literature (with some modifications) while
the last example represents a practical problem. 
Particularly, we demonstrate
the exponential rate of the convergence of the error and its estimator with respect
to the number of degrees of freedom ($\DoF$), i.e.,
$  |J(\u) - J(\uh)| \approx |\etaI (\uh,\zh)| =  O(\exp(-\DoF^{1/3}));$
cf., the theoretical results in \cite[Theorem~4.63]{Schwab-book},
\cite[Theorem~3.2]{Babuska-Suri90} and the computational results in \cite{Demko02,SolDemko04}.
Therefore, we plot the error  convergence in log-linear graphs.

\subsection{Semilinear problem}

This example has only a mild nonlinearity in the reaction term; whereas, the diffusion is
linear. 
We compare the error estimator resulting from the original differentiation of the
primal problem \eqref{i3D} and the proposed linearization in \eqref{DPh}.
Moreover, since the
solutions of the primal and dual problems do not suffer from the lack of regularity,
we demonstrate the superiority of the higher order approximations.

Let $\Om:=(0,1)^2$. In virtue of \cite[Example~35]{BeckerALL_CAMWA22},  
we consider the semilinear problem (in the weak form)
\begin{align}
  \label{TT1}
 \mbox{find }  u\in H^1_0(\Om):\quad
  \Lsp{\nabla u}{\nabla v} + \Lsp{u^3}{v} = \Lsp{\chi_{\Om_1}}{\pd v/\pd {x_1}},
\end{align}
where $\chi_{\Om_1}$ is the characteristic function of
$\Om_1:=\{x\in\Om;\, x_1+x_2\leq 0.5\}$. The quantity of interest is given
by
\begin{align}
  \label{TT2}
  J(u) =  -\Lsp{\chi_{\Om_1}}{\pd_{x_2} u},
\end{align}
where  $\chi_{\Om_2}$ is the characteristic function of  $\Om_2:=\{x\in\Om;\, x_1+x_2\geq 1.5\}$. The reference value obtained by
an ``over-kill'' computation (more than 120\,000 $\DoF$) is
$J(u)= 1.58495180882\cdot10^{-3}\pm10^{-14}$.

Following \eqref{i3D}, we define the weak adjoint problem using the differentiation of \eqref{TT1}.
Let $u\in H^1_0(\Om)$,
\begin{align}
  \label{TT3}
   \mbox{find }  z\in H^1_0(\Om):\quad
  \Lsp{\nabla v}{\nabla z} + \Lsp{3u^2 z}{v} = -\Lsp{\chi_{\Om_2}}{\pd v/\pd{x_2}} \qquad
    \forall v\in H^1_0(\Om).
\end{align}
On the other hand, the proposed linearization in \eqref{DPh} reads the following
adjoint problem. Let $u\in H^1_0(\Om)$,
\begin{align}
  \label{TT4}
   \mbox{find }  z\in H^1_0(\Om):\quad
  \Lsp{\nabla v}{\nabla z} + \Lsp{u^2 z}{v} = -\Lsp{\chi_{\Om_2}}{\pd v/\pd{x_2} } \qquad
    \forall v\in H^1_0(\Om).
\end{align}

\begin{figure} [t]
  \begin{center}
    \includegraphics[width=0.46\textwidth]{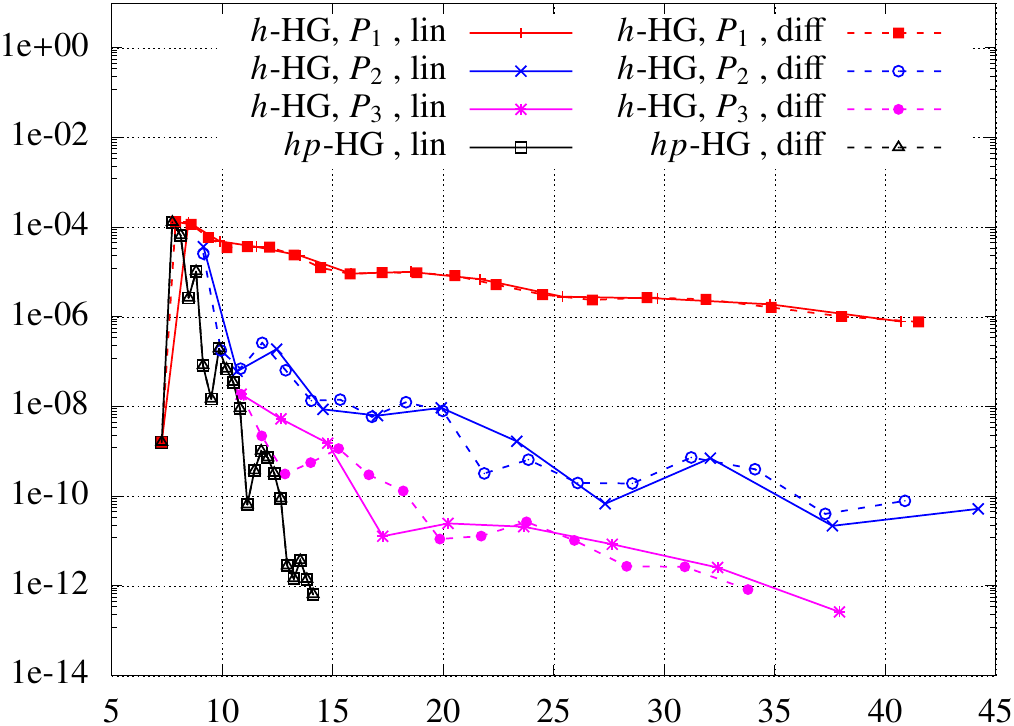}
    \hspace{0.02\textwidth}
    \includegraphics[width=0.46\textwidth]{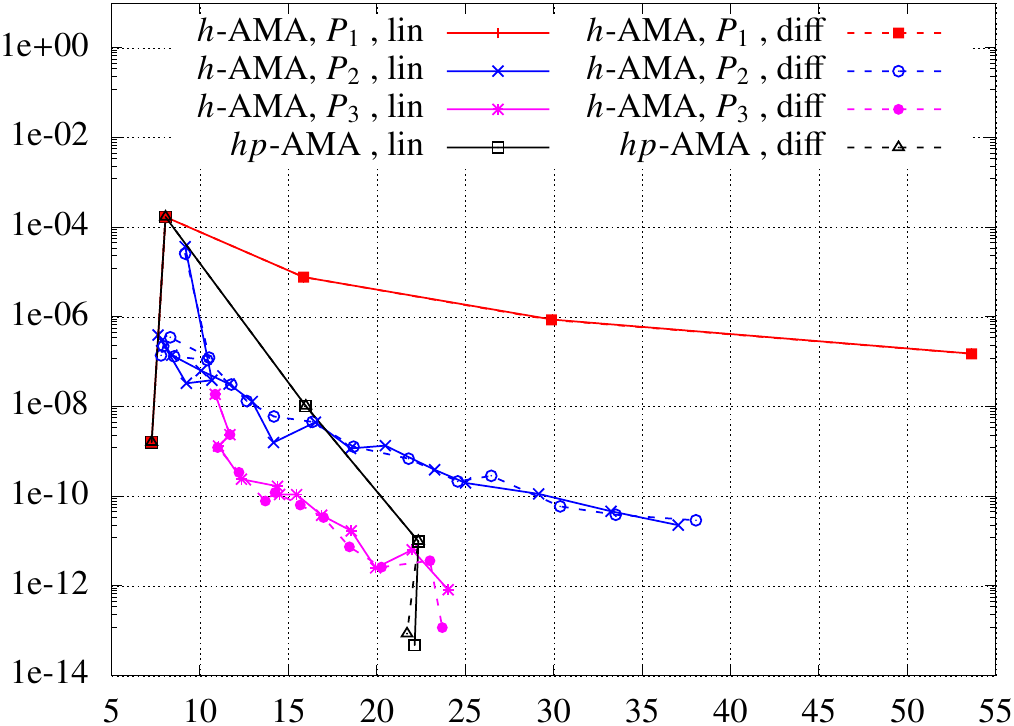}
  \end{center}
  \caption{Semilinear problem \eqref{TT1} -- \eqref{TT2}, comparison of the
    convergence of the error estimator $|\etaI|$ with respect to $\DoF^{1/3}$ for dual problems
    based on differentiation \eqref{TT3} (``diff'', dashed lines)
    and linearization  \eqref{TT4} (``lin'', full lines)
    for different adaptive techniques, {\HGh} and {\HGhp} 
    (left) and {\AMAh} and {\AMAhp} (right).}
  \label{fig:triangs_lin}
\end{figure}

In order to compare error estimates resulting from \eqref{TT3} and \eqref{TT4},
we employ the 4 adaptive techniques,
{\AMAh}, {\AMAhp}, {\HGh} and {\HGhp}, introduced in Sections~\ref{sec:AMA}--\ref{sec:HG}.
Figure~\ref{fig:triangs_lin} shows the convergence of the error estimator $|\etaI|$
for all adaptive techniques depending on the definitions of the adjoint problems
either by \eqref{TT3} or by \eqref{TT4}.
We observe very similar convergence
of error estimates using both adjoint problems for all tested adaptive techniques.
We note that the convergence of  error $|J(u)-J(\uh)|$
shows the same similarity (these graphs are not shown here). These results justify that the
use of the adjoint problem \eqref{TT4} not based on the differentiation of the primal one
is possible.

\begin{figure} 
  \begin{center}
    \includegraphics[width=0.46\textwidth]{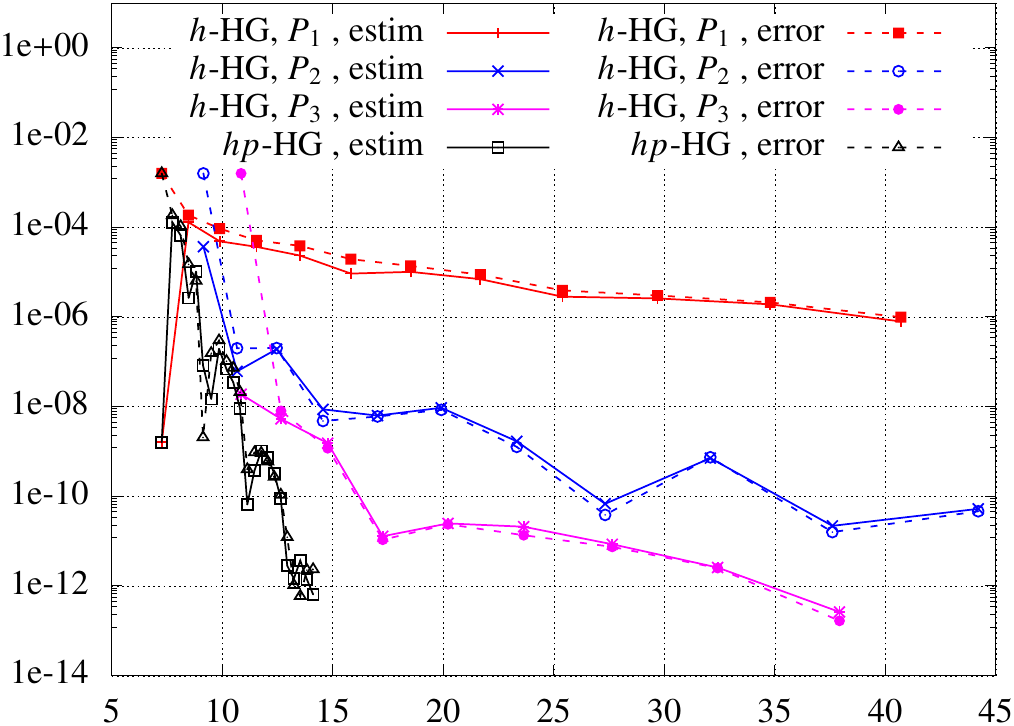}
    \hspace{0.02\textwidth}
    \includegraphics[width=0.46\textwidth]{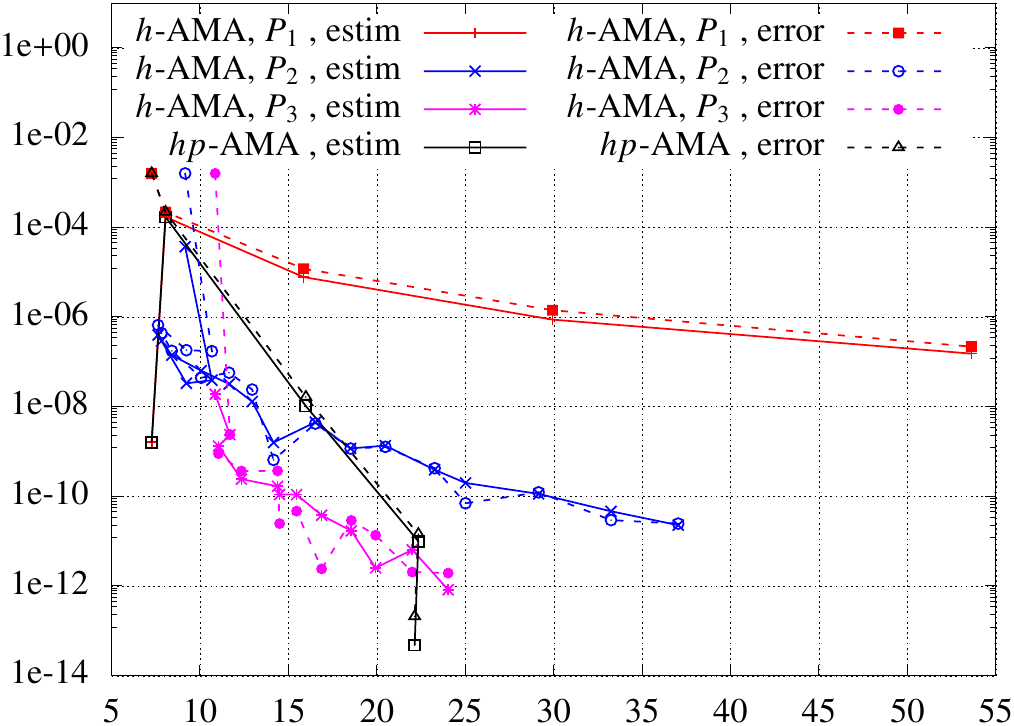}
  \end{center}
  \caption{Semilinear problem \eqref{TT1} -- \eqref{TT2}, convergence of the error $|J(u)-J(\uh)|$ (``error'', dashed lines) and its estimator
    $|\etaI|$ (``estim'', full lines) with respect to $\DoF^{1/3}$
    for different adaptive techniques, {\HGh} and {\HGhp} 
    (left) and {\AMAh} and {\AMAhp} (right).}
  \label{fig:triangs_conv}
\end{figure}

Moreover, the accuracy of the error estimator is demonstrated
by Figure~\ref{fig:triangs_conv} where we compare the error $e_h:=|J(u)-J(\uh)|$ with its estimator
$|\etaI|$
for all adaptive techniques.
\refA{These values are also given in Tables~\ref{tab:triangsHG} and \ref{tab:triangsAMA},
  where we also show  the effectivity index $\ieff:=|\etaI|/e_h$ and the computational
  time in seconds.
  The effectivity indexes are not as close to 1 as would be expected, which is caused by
  the higher-order reconstruction operator $\KrrO$ used in \eqref{CE5}. Hence, the development
  of a more accurate reconstruction working on anisotropic $hp$-meshes is still an open problem.
However, Figure~\ref{fig:triangs_conv} shows a tight approximation of the error which is not the case for the following examples}
where the nonlinearities are much
stronger. Furthermore, due to the regularity of the exact solution, the higher-order
approximations are superior to the low-order methods.
The {\HGhp} method achieved the prescribed tolerance
using fewer $\DoF$ than {\AMAhp} but it required many more adaptive cycles; hence,
the corresponding computational times are comparable.
We note that the linearization and algebraic errors ($\approx \etaL+\etaA+\etaJ$)
in this example
is negligible due to the weak nonlinearity and therefore they are not treated.

\begin{table}
  \refA{\footnotesize
   \setlength{\tabcolsep}{2.5pt}
 \begin{tabular}{|cr|ccc|r|}
 \hline
 \multicolumn{6}{|c|}{$h$-HG,    $P_1$    } \\
 \hline
  $\ell$ & $\DoF$ & 
  $ |e_h |$  &  $|\etaI|$ & $\ieff$ &time \\ 
 \hline
     0&   384&  1.58E-03&  1.60E-09& 0.00&   0.3\\
     1&   609&  1.87E-04&  1.29E-04& 0.69&   0.7\\
     2&   969&  9.27E-05&  4.86E-05& 0.52&   1.0\\
     3&  1545&  5.08E-05&  3.65E-05& 0.72&   1.6\\
     4&  2472&  3.83E-05&  2.34E-05& 0.61&   2.5\\
     5&  3966&  1.93E-05&  9.25E-06& 0.48&   4.6\\
     6&  6369&  1.36E-05&  1.01E-05& 0.74&   9.5\\
     7& 10194&  8.72E-06&  6.93E-06& 0.80&  20.2\\
     8& 16377&  3.89E-06&  2.81E-06& 0.72&  44.0\\
     9& 26196&  2.99E-06&  2.59E-06& 0.87&  89.1\\
    10& 42135&  2.10E-06&  1.90E-06& 0.91& 187.9\\
    11& 67461&  9.86E-07&  7.89E-07& 0.80& 356.2\\
 \hline
 \end{tabular}

   \setlength{\tabcolsep}{2.5pt}
 \begin{tabular}{|cr|ccc|r|}
 \hline
 \multicolumn{6}{|c|}{$h$-HG,    $P_2$    } \\
 \hline
  $\ell$ & $\DoF$ & 
  $ |e_h |$  &  $|\etaI|$ & $\ieff$ &time \\ 
 \hline
     0&   768&  1.58E-03&  3.66E-05& 0.02&   0.4\\
     1&  1218&  2.00E-07&  6.09E-08& 0.30&   0.9\\
     2&  1938&  2.02E-07&  1.91E-07& 0.94&   1.5\\
     3&  3090&  4.78E-09&  8.69E-09& 1.82&   2.6\\
     4&  4944&  6.00E-09&  6.31E-09& 1.05&   5.1\\
     5&  7896&  8.45E-09&  9.41E-09& 1.11&   9.9\\
     6& 12702&  1.26E-09&  1.69E-09& 1.34&  19.9\\
     7& 20406&  3.84E-11&  6.82E-11& 1.77&  40.1\\
     8& 33042&  7.41E-10&  7.11E-10& 0.96&  87.2\\
     9& 53202&  1.59E-11&  2.18E-11& 1.38& 168.8\\
    10& 86376&  4.64E-11&  5.26E-11& 1.13& 417.1\\
 \hline
 \hline
 \multicolumn{6}{|c|}{$h$-HG,    $P_3$    } \\
 \hline
  $\ell$ & $\DoF$ & 
  $ |e_h |$  &  $|\etaI|$ & $\ieff$ &time \\ 
 \hline
     0&  1280&  1.58E-03&  1.88E-08& 0.00&   0.5\\
     1&  2030&  8.02E-09&  5.33E-09& 0.66&   1.7\\
     2&  3230&  1.18E-09&  1.52E-09& 1.29&   3.3\\
     3&  5150&  1.07E-11&  1.28E-11& 1.20&   6.3\\
     4&  8270&  2.38E-11&  2.48E-11& 1.04&  14.4\\
     5& 13220&  1.36E-11&  2.09E-11& 1.54&  26.8\\
     6& 21140&  7.47E-12&  8.51E-12& 1.14&  54.0\\
     7& 34070&  2.50E-12&  2.58E-12& 1.03& 100.1\\
     8& 54590&  1.67E-13&  2.64E-13& 1.58& 177.3\\
 \hline
 \end{tabular}

   \setlength{\tabcolsep}{2.5pt}
 \begin{tabular}{|cr|ccc|r|}
 \hline
 \multicolumn{6}{|c|}{$hp$-HG             } \\
 \hline
  $\ell$ & $\DoF$ & 
  $ |e_h |$  &  $|\etaI|$ & $\ieff$ &time \\ 
 \hline
     0&   384&  1.58E-03&  1.60E-09& 0.00&   0.3\\
     1&   465&  1.83E-04&  1.29E-04& 0.71&   0.5\\
     2&   537&  1.04E-04&  6.55E-05& 0.63&   0.6\\
     3&   612&  1.48E-05&  2.64E-06& 0.18&   0.8\\
     4&   687&  6.31E-06&  1.03E-05& 1.63&   0.9\\
     5&   764&  2.05E-09&  8.19E-08&39.94&   1.1\\
     6&   860&  1.55E-07&  1.47E-08& 0.10&   1.3\\
     7&   960&  2.95E-07&  1.94E-07& 0.66&   1.5\\
     8&  1060&  1.03E-07&  6.85E-08& 0.67&   1.8\\
     9&  1160&  7.23E-08&  3.48E-08& 0.48&   2.2\\
    10&  1265&  2.10E-08&  9.10E-09& 0.43&   2.6\\
    11&  1384&  3.95E-10&  6.54E-11& 0.17&   3.0\\
    12&  1507&  9.39E-10&  3.72E-10& 0.40&   3.6\\
    13&  1631&  9.21E-10&  1.02E-09& 1.10&   4.2\\
    14&  1757&  5.90E-10&  7.20E-10& 1.22&   4.8\\
    15&  1888&  2.69E-10&  3.26E-10& 1.21&   5.6\\
    16&  2028&  1.07E-10&  8.83E-11& 0.83&   6.5\\
    17&  2172&  1.22E-11&  2.82E-12& 0.23&   7.5\\
    18&  2325&  1.04E-12&  1.46E-12& 1.40&   8.7\\
    19&  2482&  5.93E-13&  3.68E-12& 6.20&  10.3\\
    20&  2648&  2.25E-12&  1.41E-12& 0.63&  12.0\\
    21&  2816&  2.35E-12&  6.48E-13& 0.28&  14.0\\
 \hline
 \end{tabular}

  }
  \caption{Semilinear problem \eqref{TT1} -- \eqref{TT2},  error $e_h=|J(u)-J(\uh)|$,
    its estimate $|\etaI|$, effectivity index $\ieff$ and the computational time in seconds
  for {\HGh} and {\HGhp} adaptive methods.}
  \label{tab:triangsHG}
\end{table}

\begin{table}
  \refA{\footnotesize
       \setlength{\tabcolsep}{2.5pt}
 \begin{tabular}{|cr|ccc|r|}
 \hline
 \multicolumn{6}{|c|}{$h$-AMA,   $P_1$    } \\
 \hline
  $\ell$ & $\DoF$ & 
  $ |e_h |$  &  $|\etaI|$ & $\ieff$ &time \\ 
 \hline
     0&   384&  1.58E-03&  1.60E-09& 0.00&   0.4\\
     1&   522&  2.17E-04&  1.69E-04& 0.78&   0.6\\
     2&  3987&  1.19E-05&  7.71E-06& 0.65&   3.3\\
     3& 26832&  1.40E-06&  8.67E-07& 0.62&  37.8\\
     4&154239&  2.16E-07&  1.51E-07& 0.70& 704.7\\
 \hline
 \end{tabular}

       \setlength{\tabcolsep}{2.5pt}
 \begin{tabular}{|cr|ccc|r|}
 \hline
 \multicolumn{6}{|c|}{$h$-AMA,   $P_2$    } \\
 \hline
  $\ell$ & $\DoF$ & 
  $ |e_h |$  &  $|\etaI|$ & $\ieff$ &time \\ 
 \hline
     0&   768&  1.58E-03&  3.66E-05& 0.02&   0.4\\
     1&  1218&  1.72E-07&  3.89E-08& 0.23&   0.9\\
     2&   786&  1.81E-07&  3.31E-08& 0.18&   1.1\\
     3&   444&  6.53E-07&  3.94E-07& 0.60&   1.3\\
     4&   480&  4.32E-07&  2.89E-07& 0.67&   1.4\\
     5&   594&  1.74E-07&  1.39E-07& 0.80&   1.6\\
     6&  1020&  4.38E-08&  6.32E-08& 1.44&   1.9\\
     7&  1578&  5.70E-08&  3.18E-08& 0.56&   2.4\\
     8&  2172&  2.39E-08&  1.29E-08& 0.54&   3.1\\
     9&  2838&  6.50E-10&  1.58E-09& 2.42&   4.1\\
    10&  4512&  4.16E-09&  4.47E-09& 1.07&   5.9\\
    11&  6354&  1.15E-09&  1.17E-09& 1.01&   8.8\\
    12&  8598&  1.27E-09&  1.34E-09& 1.06&  13.3\\
    13& 12606&  4.24E-10&  3.89E-10& 0.92&  20.7\\
    14& 15654&  6.99E-11&  1.98E-10& 2.84&  30.8\\
    15& 24780&  1.23E-10&  1.13E-10& 0.92&  50.2\\
    16& 36750&  2.94E-11&  4.64E-11& 1.58&  84.5\\
    17& 50820&  2.46E-11&  2.27E-11& 0.92& 144.6\\
 \hline
 \end{tabular}

       \setlength{\tabcolsep}{2.5pt}
 \begin{tabular}{|cr|ccc|r|}
 \hline
 \multicolumn{6}{|c|}{$h$-AMA,   $P_3$    } \\
 \hline
  $\ell$ & $\DoF$ & 
  $ |e_h |$  &  $|\etaI|$ & $\ieff$ &time \\ 
 \hline
     0&  1280&  1.58E-03&  1.88E-08& 0.00&   0.6\\
     1&  1600&  2.38E-09&  2.37E-09& 0.99&   1.3\\
     2&  1350&  8.90E-10&  1.32E-09& 1.48&   1.8\\
     3&  1880&  3.65E-10&  2.41E-10& 0.66&   2.6\\
     4&  2980&  3.69E-10&  1.68E-10& 0.45&   3.8\\
     5&  3050&  2.46E-11&  1.10E-10& 4.49&   5.1\\
     6&  3700&  4.67E-11&  1.09E-10& 2.34&   6.7\\
     7&  4810&  2.41E-12&  3.76E-11&15.63&   9.0\\
     8&  6390&  2.89E-11&  1.69E-11& 0.58&  12.2\\
     9&  7920&  1.36E-11&  2.54E-12& 0.19&  16.4\\
    10& 10650&  2.05E-12&  6.41E-12& 3.13&  24.5\\
    11& 13890&  1.93E-12&  8.24E-13& 0.43&  33.1\\
 \hline
 \hline
 \multicolumn{6}{|c|}{$hp$-AMA            } \\
 \hline
  $\ell$ & $\DoF$ & 
  $ |e_h |$  &  $|\etaI|$ & $\ieff$ &time \\ 
 \hline
     0&   384&  1.58E-03&  1.60E-09& 0.00&   0.4\\
     1&   522&  2.17E-04&  1.69E-04& 0.78&   0.8\\
     2&  4086&  1.63E-08&  1.03E-08& 0.63&   3.1\\
     3& 11162&  1.43E-11&  9.93E-12& 0.69&  10.0\\
     4& 10850&  2.08E-13&  4.82E-14& 0.23&  17.5\\
 \hline
 \end{tabular}

  }
  \caption{Semilinear problem \eqref{TT1} -- \eqref{TT2},  error $e_h=|J(u)-J(\uh)|$,
    its estimate $|\etaI|$, effectivity index $\ieff$ and the computational time in seconds
  for {\AMAh} and {\AMAhp} adaptive methods.}
  \label{tab:triangsAMA}
\end{table}

\subsection{Quasilinear elliptic problem on L-shaped domain}
\begin{figure}
  \begin{center}
    \includegraphics[height=0.21\textwidth]{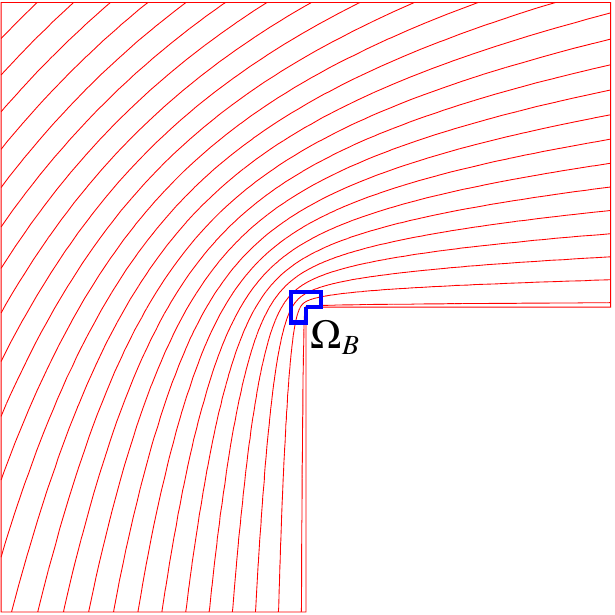}
    \hspace{1mm}
    \includegraphics[height=0.21\textwidth]{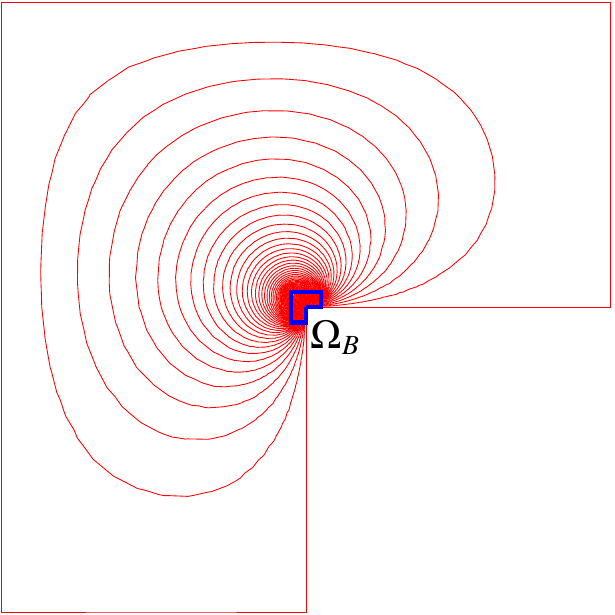}
     \hspace{1mm}
     \includegraphics[height=0.21\textwidth]{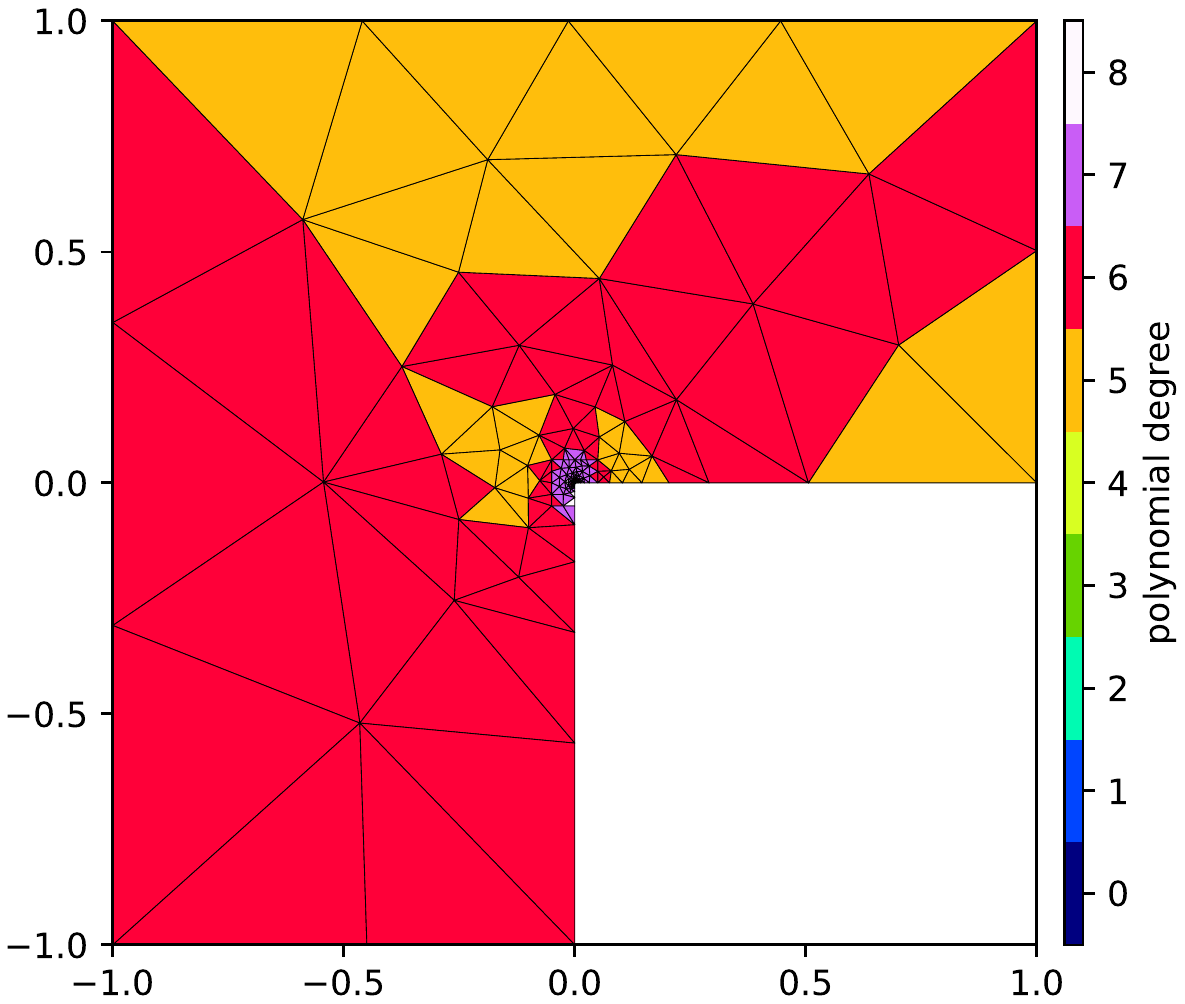}
    \hspace{1mm}
    \includegraphics[height=0.21\textwidth]{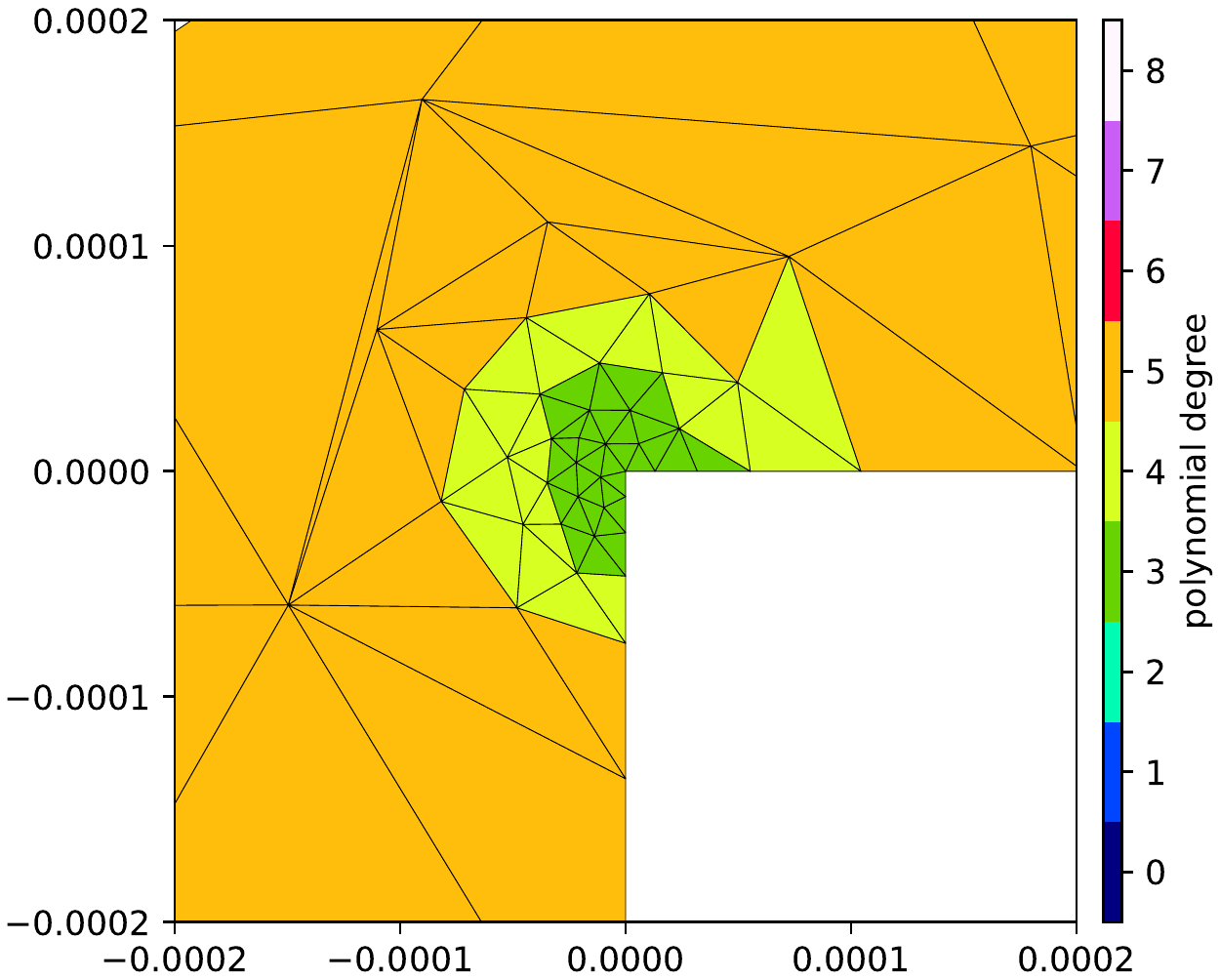}
  \end{center}
  \caption{Quasilinear elliptic problem \eqref{Q1}--\eqref{Q2}:
    primal solution (first), adjoint solution (second), and the final $hp$-mesh (third)
    and its 5000x zoom  near the interior corner (last).}
  \label{fig:quasiHH_init}
\end{figure}

Similarly as in \cite{houston-robson-suli}, we consider a
quasilinear elliptic problem on L-shaped domain
$\Om:=(-1, 1)^2 \setminus [0, 1)\times (-1, 0]$  
\begin{align}
  \label{Q1}
  -\nabla \cdot( \mu(|\nabla u|) \nabla u) &=  f\qquad\hspace{7pt}\mbox{ in }\Om, 
\end{align}
where $\mu(|\nabla u|) = 1 + \exp(-|\nabla u|^2)$ is a nonlinear diffusion.
We prescribe the Dirichlet boundary condition on the boundary $\pd\Om$ and 
the function $f$ such that the
exact solution is $u = r^{2/3}\sin(2\phi/3)$ with $(r,\phi)$ being the polar coordinates.
The target functional represents the total energy (cf.~\eqref{alt7} in Section~\ref{sec:alter})
in a small polygonal domain around the interior corner
$\Om_B:= \{ (x_1,x_2)\in \Om,\ |x_i| \leq 0.05,\ i=1,2 \}$; i.e.,
\begin{align}
  \label{Q2}
  J(\u) = \int\nolimits_{\Om_B} \mu(|\nabla u|) |\nabla u|^2\dx.
\end{align}
Since the exact solution is know, the reference value
$J(\u)= 1.721609238808\cdot10^{-2}$ has been computed using a numerical quadrature.
Figure~\ref{fig:quasiHH_init}, left, shows the primal and adjoint solutions corresponding to
\eqref{Q1}--\eqref{Q2} together with 
the domain of interest $\Om_B$.
Due to the interior angle the primal and adjoint solutions have a singularity.

\begin{figure}
  \begin{center}
    \includegraphics[height=0.29\textwidth]{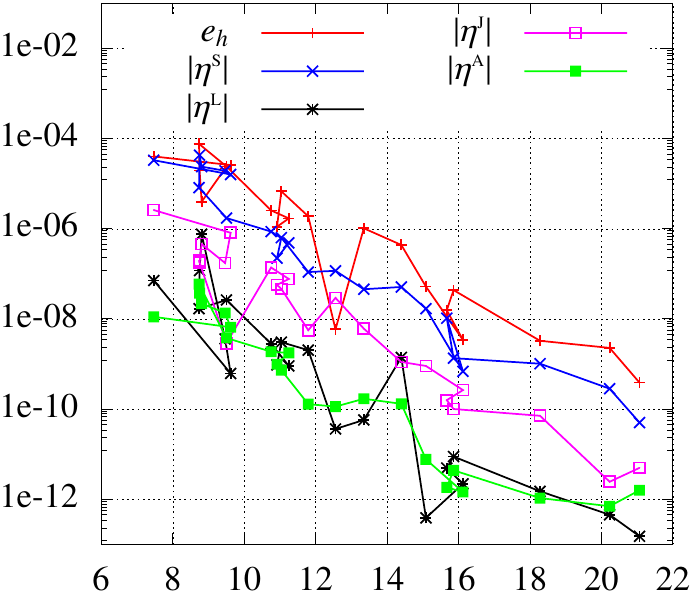}
    \hspace{8mm}
    \includegraphics[height=0.29\textwidth]{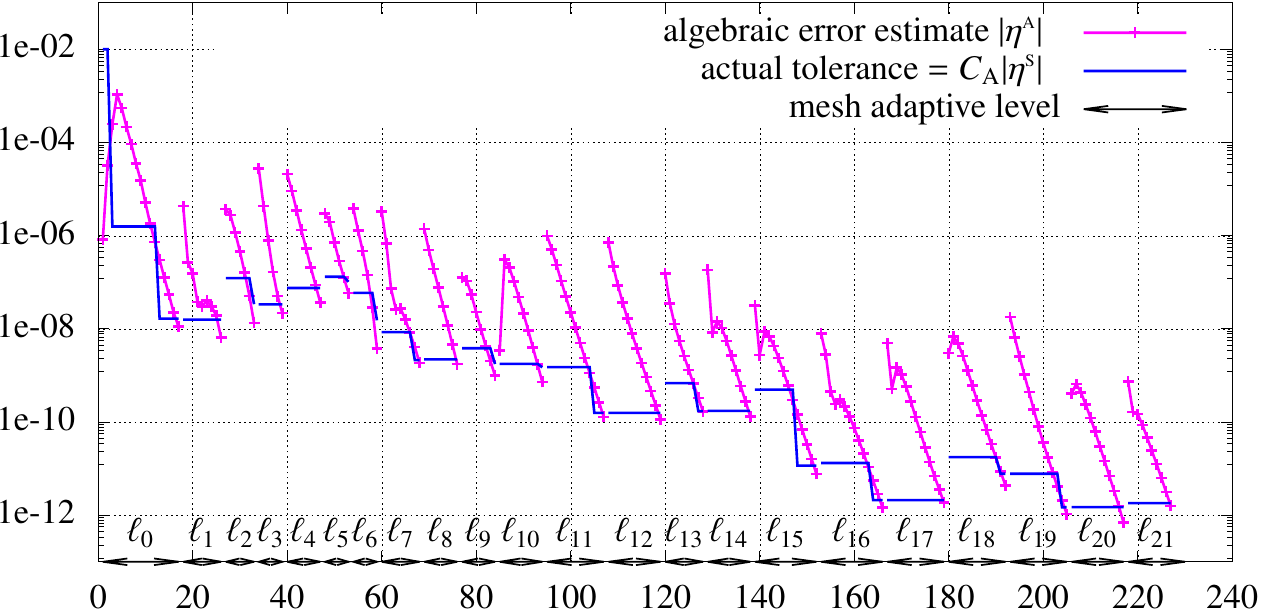}
  \end{center}
  \caption{Quasilinear elliptic problem \eqref{Q1}--\eqref{Q2},
    convergence of $e_h=|J(\u)-J(\uh)|$ and its estimators $|\etaS|$, $|\etaL|$, $|\etaJ|$
    and $|\etaA|$
    with respect to $\DoF^{1/3}$
    (left) and the convergence of the nonlinear
    iterative solver for adaptive levels $\ell=0,\dots 21$.}
  \label{fig:quasiHH_conv}
\end{figure}

We solve this problem only by the {\AMAhp} adaptive technique.
Figure~\ref{fig:quasiHH_conv}, left,  shows
the convergence of the error $e_h:=|J(\u)-J(\uh)|$ and the various parts of the
error estimator: linearization $|\etaL|$ and $|\etaJ|$, residual $|\etaS|$, and algebraic $|\etaA|$
with respect to $\DoF^{1/3}$, cf.~\eqref{CE10a};
each node corresponds to one level of mesh adaptation.
\refA{Moreover, Table~\ref{tab:quasiHH} shows the vales of $e_h$, $|\etaI|$, the effectivity index
  $\ieff:=|\etaI|/e_h$ and the computational time in seconds.}
As we do
not have an upper bound of the error, it is underestimated at most by a factor of 10.
The dominating term is the residual estimator $|\etaS|$.
The convergence of all estimators is not monotone, which is a typical behaviour
for anisotropic adaptation.
Nevertheless, an exponential rate of convergence is observed.
The prescribed tolerance $\TOL=10^{-10}$ is achieved after 21 levels of mesh adaptation.

Moreover, Figure~\ref{fig:quasiHH_conv}, right, shows
the convergence of the nonlinear solver 
level $\ell=0,1,\dots, 21$ of mesh adaptation. 
Particularly, we plot both sides of the stopping criterion
\eqref{CE15}, i.e., the estimate of the algebraic error $|\etaA(\uhlk,\zhlk)|$ and the
adaptively chosen tolerance $\CA |\etaS|$ (with $\CA=0.1$).
Each node corresponds to one nonlinear iteration $k=0,1,\dots$, but on different adaptive level
$\ell=0,1,\dots$, in general.
The black left-right arrows indicates one mesh adaptive loop $\ell$.
We observe that the algebraic error tolerance  $\CA  |\etaS| $ only requires recalculation at most once within
one mesh adaptive level (for $\ell=2,6,7,11,13,15,16,19$). Further, this tolerance is
step by step decreasing (but not, in general, monotonic) when the total error estimates
$|\etaI|=|\etaS+\etaA+\etaL+\etaJ|$ is approaching to the tolerance $\TOL$.

\begin{table}
  \refA{\footnotesize
       \setlength{\tabcolsep}{2.5pt}
 \begin{tabular}{|cr|ccc|r|}
 \hline
 \multicolumn{6}{|c|}{$hp$-AMA            } \\
 \hline
  $\ell$ & $\DoF$ & 
  $ |e_h |$  &  $|\etaI|$ & $\ieff$ &time \\ 
 \hline
     0&   417&  3.97E-05&  3.26E-05& 0.82&   0.9\\
     1&   891&  2.56E-05&  1.58E-05& 0.62&   1.9\\
     2&   849&  2.14E-05&  1.90E-05& 0.89&   3.5\\
     3&   684&  3.92E-06&  2.35E-05& 5.99&   4.4\\
     4&   671&  1.93E-05&  4.24E-05& 2.19&   5.4\\
     5&   668&  7.48E-05&  8.06E-06& 0.11&   6.7\\
     6&   859&  2.37E-05&  1.72E-06& 0.07&   8.7\\
     7&  1244&  2.53E-06&  8.58E-07& 0.34&  12.0\\
 \hline
 \end{tabular}
 \hspace{1mm}
 \setlength{\tabcolsep}{2.5pt}
 \begin{tabular}{|cr|ccc|r|}
 \hline
 \multicolumn{6}{|c|}{$hp$-AMA            } \\
 \hline
  $\ell$ & $\DoF$ & 
  $ |e_h |$  &  $|\etaI|$ & $\ieff$ &time \\ 
 \hline
     8&  1425&  1.68E-06&  4.84E-07& 0.29&  14.5\\
     9&  1303&  1.09E-06&  2.22E-07& 0.20&  19.0\\
    10&  1346&  6.70E-06&  6.29E-07& 0.09&  23.2\\
    11&  1642&  1.86E-06&  1.09E-07& 0.06&  28.7\\
    12&  1982&  5.86E-09&  1.16E-07&19.82&  33.8\\
    13&  2379&  1.02E-06&  4.54E-08& 0.04&  44.1\\
    14&  2987&  4.32E-07&  5.07E-08& 0.12&  55.4\\
    15&  3438&  5.15E-08&  1.67E-08& 0.32&  79.3\\
 \hline
 \end{tabular}
 \hspace{1mm}
 \setlength{\tabcolsep}{2.5pt}
 \begin{tabular}{|cr|ccc|r|}
 \hline
 \multicolumn{6}{|c|}{$hp$-AMA            } \\
 \hline
  $\ell$ & $\DoF$ & 
  $ |e_h |$  &  $|\etaI|$ & $\ieff$ &time \\ 
 \hline
    16&  4193&  3.40E-09&  6.84E-10& 0.20& 116.2\\
    17&  3851&  1.55E-08&  1.03E-08& 0.67& 144.0\\
    18&  3993&  4.34E-08&  1.35E-09& 0.03& 186.0\\
    19&  6111&  3.25E-09&  1.01E-09& 0.31& 252.0\\
    20&  8278&  2.26E-09&  2.83E-10& 0.13& 317.7\\
    21&  9348&  3.85E-10&  5.01E-11& 0.13& 412.2\\
 \hline
 \end{tabular}

  }
  \caption{Quasilinear elliptic problem \eqref{Q1}--\eqref{Q2}, error $e_h=|J(u)-J(\uh)|$,
    its estimate $|\etaI|$, effectivity index $\ieff$ and the computational time in seconds
    obtained by {\AMAhp} adaptive method for adaptive levels $\ell=0,\dots 21$.}
  \label{tab:quasiHH}
\end{table}

The resulting $hp$-grid with detail near the interior corner is shown in
Figure~\ref{fig:quasiHH_init}, right.
The large elements with high polynomial degrees are outside
of the singularity; whereas, small elements with a low polynomial degree ($p=3$) are generated
near the interior corner.

\subsection{Convective-dominated problem with the Carreau-law diffusion}

We consider the convection-diffusion problem in the form 
\begin{align}
  \label{car1}
  -\nabla\cdot(\mu(|\nabla u|) \nabla u) + \bkb \cdot \nabla u = 0
  \qquad \mbox{ in } \Om:= (0,2)\times(0,1),
\end{align}
where $\bkb = (x_2 , (1-x_1)^2 )$ is the prescribed velocity field and the nonlinear diffusion
is given by the Carreau law for a non-Newtonian fluid
(\cite{BarrettLiu_NM93,BerroneSuli_IMA08,CongreveALL_IMA13})
\begin{align}
  \label{car2}
  \mu(|\nabla u|) = \ve(k_\infty+(k_0 - k_\infty))
  \left(1 + \lambda |\nabla u|^2\right)^{(\theta-2)/\theta},
\end{align}
where $\ve>0$, $\lambda > 0$, $1< \theta \leq 2$ and $0 < k_\infty < k_0$ .
We prescribe the homogeneous Neumann data at the outflow part
$\gom_N := \{2\}\times(0,1) \cup (0,2)\times \{1\} $ 
and the discontinuous Dirichlet data 
\begin{equation}
  \label{car3}
  u = \left\{ 
  \begin{array}{ll}
  1 & x_1 \in (\frac18, \frac12),\ x_2 = 0 \\[2pt]
  2 & x_1 \in (\frac12, \frac34),\ x_2 = 0 \\[2pt]
 0 &  \mbox{ elsewhere on } \gom_D:=\gom\setminus\gom_N.
  \end{array}
  \right.
\end{equation}
We consider the values $\ve=10^{-4}$, $\lambda = 1$, $\theta=1.2$, $k_\infty = 1$ and
$k_0=2$. %
The discontinuity of the boundary conditions leads to the presence of three interior layers
which propagates through the computational domain and which are smeared due to the presence of
diffusion.

The quantity of interest is given by the integral 
$  J(u) = \int_{\Gamma_B} u\dS$, where
$\Gamma_B:=\{ (x_1,x_2), \ 1.8<x_1<2,\ x_2 = 1\} \cup \{ (x_1,x_2),\ x_1=1,\ 0.8<x_2<1\}$
is a part of the Neumann boundary $\gom_N$.
The reference value obtain by computations obtained on a strongly refined grid is
$0.280006172$.
We solved this problem with the {\AMAhp} method where the error tolerance is $\TOL=10^{-10}$.

\begin{figure} 
  \begin{center}
    \includegraphics[height=0.29\textwidth]{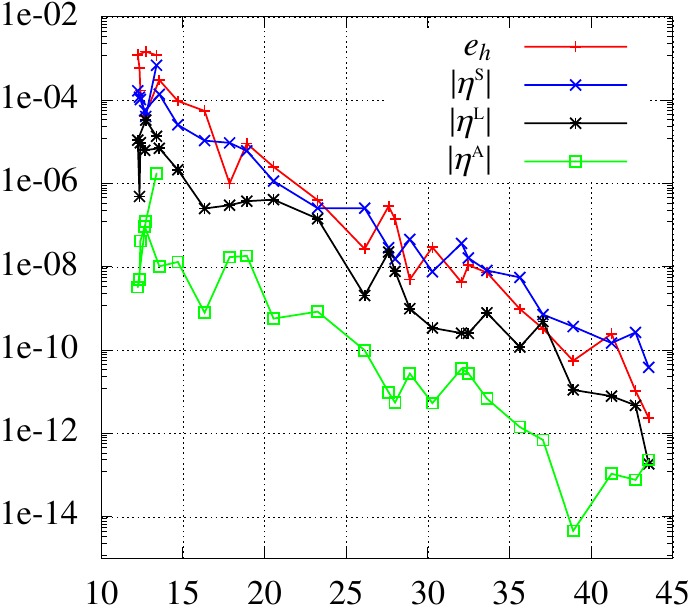}
    \hspace{0.04\textwidth}
    \includegraphics[height=0.29\textwidth]{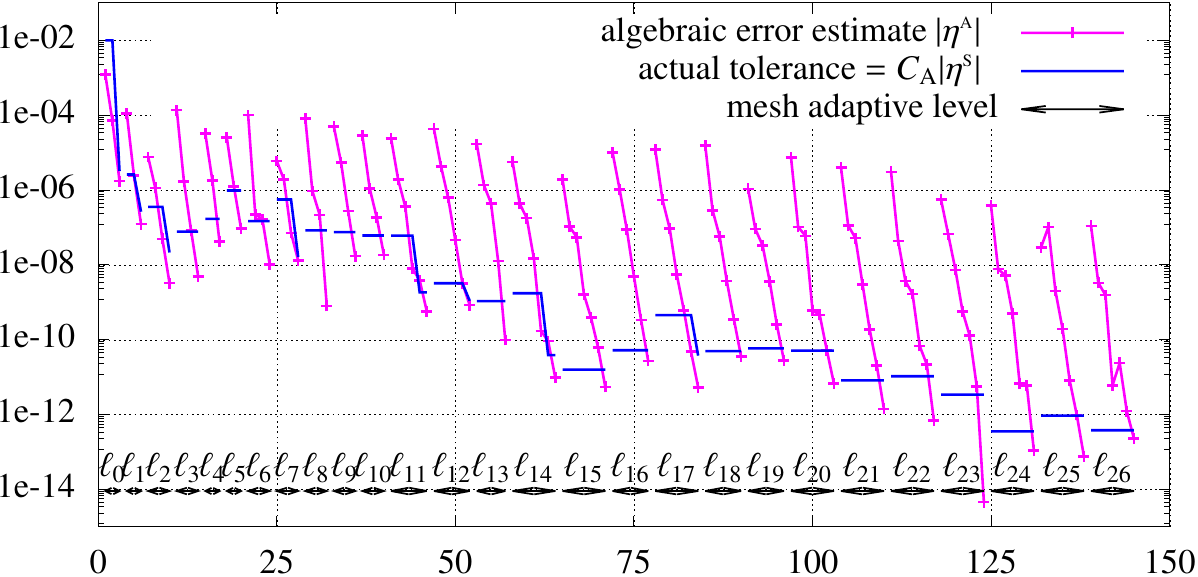}
  \end{center}
  \caption{Convection with the Carreau law diffusion \eqref{car1}--\eqref{car3},
    convergence of the error $e_h=|J(u)-J(u_h)|$
    and its estimators $|\etaS|$, $|\etaL|$ and $|\etaA|$
    with respect to $\DoF^{1/3}$ (left) and the convergence of the nonlinear
    iterative solver for adaptive levels $\ell=0,1,\dots,26$.}
  \label{fig:carreaur_conv}
\end{figure}

Figure~\ref{fig:carreaur_conv}, left, shows the convergence of the error $e_h=|J(u)-\J(\uh)|$
and  the various parts of the
error estimator: linearization $|\etaL|$, residual $|\etaS|$ and algebraic $|\etaA|$
with respect to $\DoF^{1/3}$, cf.~\eqref{CE10a}; each node corresponds
to one level of mesh adaptation. Here,  $\etaJ=0$.
\refA{The values $e_h$, $\etaI$, effectivity index $\ieff$ and computational time
are shown in Table~\ref{tab:carreau}.}
We observe the exponential rate of the convergence and
a reasonable approximation of the error,
\refA{about $0.1\lesssim\ieff\lesssim 10$ due to a strong anisotropy of the meshes.
  It is obvious namely for the last
  two levels of adaption where the limits of finite
  precision arithmetic and the error in the reference
  value for the quantity of interest, obtained by a highly refined mesh approximation,
  may both play non-negligible roles.}
The dominant part of the estimator is $|\etaS|$;
however,
the role of the estimator $|\etaL|$ is larger in comparison to previous
examples. 

Moreover, Figure~\ref{fig:carreaur_conv}, right, shows the convergence of
the nonlinear solver 
for each
mesh adaptation level $\ell=0,1,\dots, 26$. 
Again, we plot both sides of the stopping criterion \eqref{CE15},
the algebraic estimate $|\etaA(\uhlk,\zhlk)|$ and the
adaptively chosen tolerance $\CA |\etaS|$ (with $\CA=0.1$).
Each node corresponds to one nonlinear iteration $k=0,1,\dots$ but on a different adaptive level
$\ell=0,1,\dots$, in general.
The black left-right arrows indicates one mesh adaptive loop $\ell$.
A step by step decrease of the error tolerance and the convergence of the nonlinear
solver is obvious.

\begin{table}
  \refA{\footnotesize
       \setlength{\tabcolsep}{2.5pt}
 \begin{tabular}{|cr|ccc|r|}
 \hline
 \multicolumn{6}{|c|}{$hp$-AMA            } \\
 \hline
  $\ell$ & $\DoF$ & 
  $ |e_h |$  &  $|\etaI|$ & $\ieff$ &time \\ 
 \hline
     0&  2400&  1.20E-03&  6.87E-04& 0.57&   2.3\\
     1&  2064&  1.46E-03&  4.04E-05& 0.03&   4.1\\
     2&  1833&  1.23E-03&  1.68E-04& 0.14&   5.8\\
     3&  1881&  5.90E-04&  1.17E-04& 0.20&   7.2\\
     4&  1908&  1.70E-04&  1.04E-04& 0.61&   8.5\\
     5&  2033&  4.23E-05&  5.61E-05& 1.33&   9.9\\
     6&  2493&  3.05E-04&  1.39E-04& 0.46&  11.8\\
     7&  3177&  9.51E-05&  2.55E-05& 0.27&  15.0\\
     8&  4360&  5.51E-05&  1.07E-05& 0.19&  18.6\\
     9&  5702&  9.94E-07&  9.38E-06& 9.44&  23.6\\
 \hline
 \end{tabular}
 \hspace{1mm}
 \setlength{\tabcolsep}{2.5pt}
 \begin{tabular}{|cr|ccc|r|}
 \hline
 \multicolumn{6}{|c|}{$hp$-AMA            } \\
 \hline
  $\ell$ & $\DoF$ & 
  $ |e_h |$  &  $|\etaI|$ & $\ieff$ &time \\ 
 \hline
    10&  6773&  9.05E-06&  6.07E-06& 0.67&  30.0\\
    11&  8684&  2.51E-06&  1.14E-06& 0.45&  42.3\\
    12& 12574&  3.97E-07&  2.55E-07& 0.64&  63.7\\
    13& 17848&  2.67E-08&  2.54E-07& 9.51&  88.4\\
    14& 21028&  2.82E-07&  2.85E-08& 0.10& 143.9\\
    15& 21948&  1.39E-07&  1.53E-08& 0.11& 187.0\\
    16& 24134&  4.96E-09&  4.55E-08& 9.18& 234.0\\
    17& 27754&  2.93E-08&  7.62E-09& 0.26& 324.1\\
    18& 32985&  4.25E-09&  3.62E-08& 8.52& 393.0\\
    19& 34260&  1.09E-08&  1.62E-08& 1.48& 471.5\\
 \hline
 \end{tabular}
 \hspace{1mm}
 \setlength{\tabcolsep}{2.5pt}
 \begin{tabular}{|cr|ccc|r|}
 \hline
 \multicolumn{6}{|c|}{$hp$-AMA            } \\
 \hline
  $\ell$ & $\DoF$ & 
  $ |e_h |$  &  $|\etaI|$ & $\ieff$ &time \\ 
 \hline
    20& 37989&  7.22E-09&  8.13E-09& 1.13& 569.4\\
    21& 45251&  9.66E-10&  5.55E-09& 5.74& 697.0\\
    22& 50881&  3.28E-10&  7.10E-10& 2.16& 853.0\\
    23& 58855&  5.53E-11&  3.68E-10& 6.65&1048.0\\
    24& 70198&  2.45E-10&  1.47E-10& 0.60&1309.4\\
    25& 77932&  1.05E-11&  2.64E-10&25.23&1646.3\\
    26& 82447&  2.36E-12&  3.83E-11&16.25&2058.4\\
 \hline
 \end{tabular}

  }
  \caption{Convection with the Carreau law diffusion \eqref{car1}--\eqref{car3},
    error $e_h=|J(u)-J(\uh)|$,
    its estimate $|\etaI|$, effectivity index $\ieff$ and the computational time in seconds
    obtained by {\AMAhp} adaptive method for adaptive levels $\ell=0,\dots 26$.}
  \label{tab:carreau}
\end{table}

\begin{figure} 
  \begin{center}
    \includegraphics[width=0.48\textwidth]{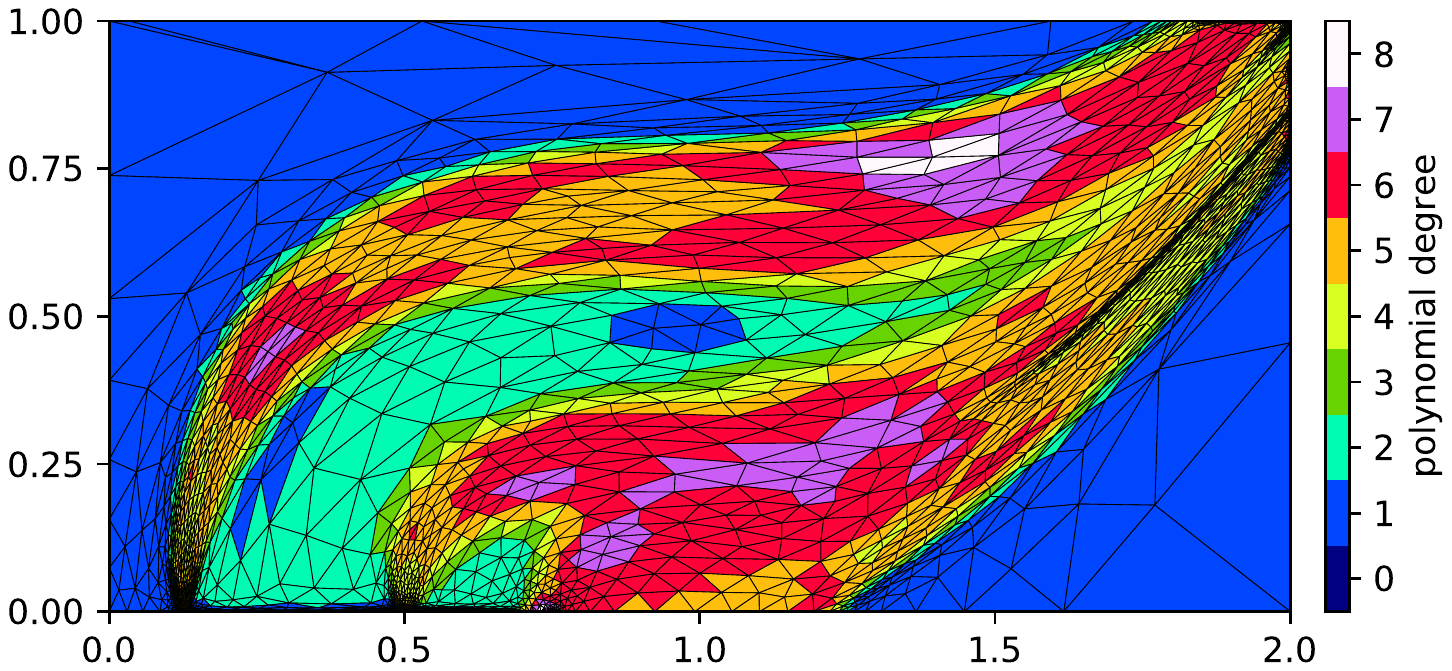}
    \hspace{0.02\textwidth}
    \includegraphics[width=0.48\textwidth]{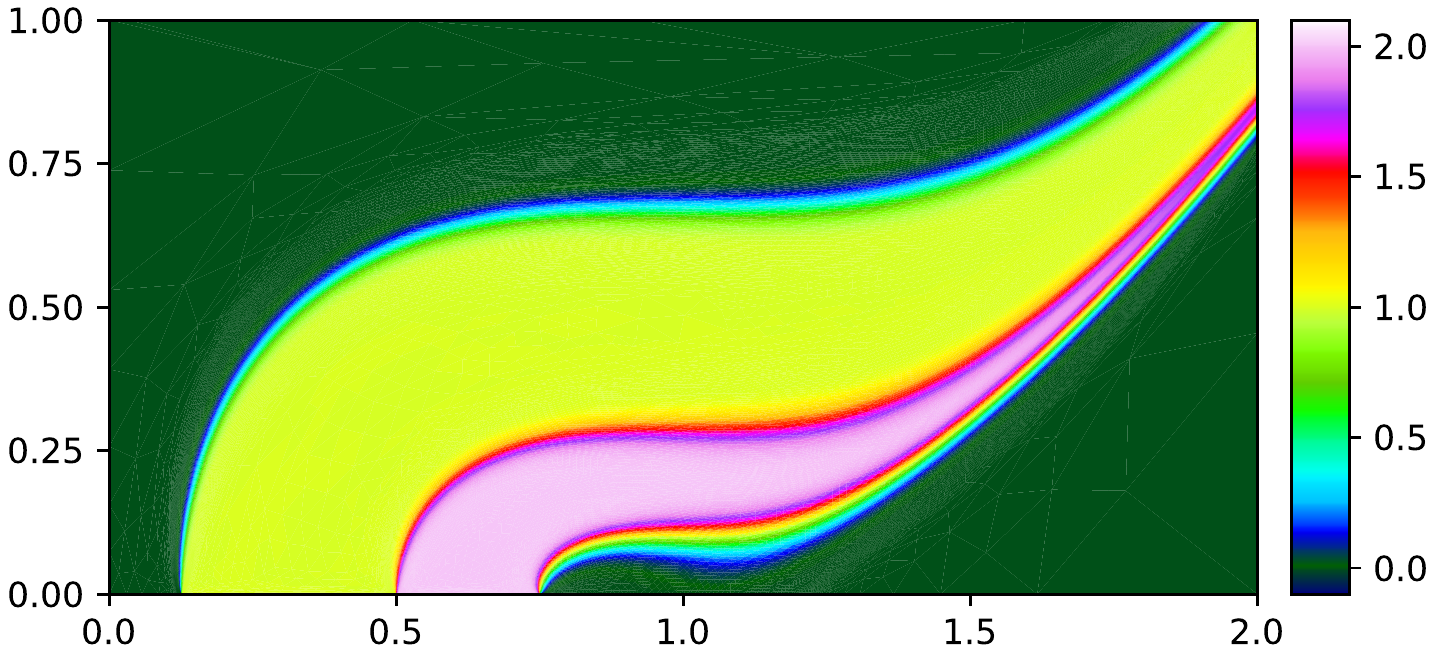}
  \end{center}
  \caption{Convection with the Carreau law diffusion \eqref{car1}--\eqref{car3},
    final $hp$-mesh (left)
    and the solution (right) obtained by {\AMAhp} method.}
  \label{fig:carreaur_meshes}
\end{figure}

Finally, Figure~\ref{fig:carreaur_meshes} shows the resulting $hp$-grid and the
corresponding solution. A mesh alignment of anisotropic elements along the interior layers
is obvious; the lowest polynomial degree and larger mesh elements are generated outside of these
layers, where the solution is almost constant.

\subsection{Magneto-static field of an alternator}
\label{sec:alter}

\begin{figure} [t]
  \begin{center}
    \includegraphics[height=0.34\textwidth]{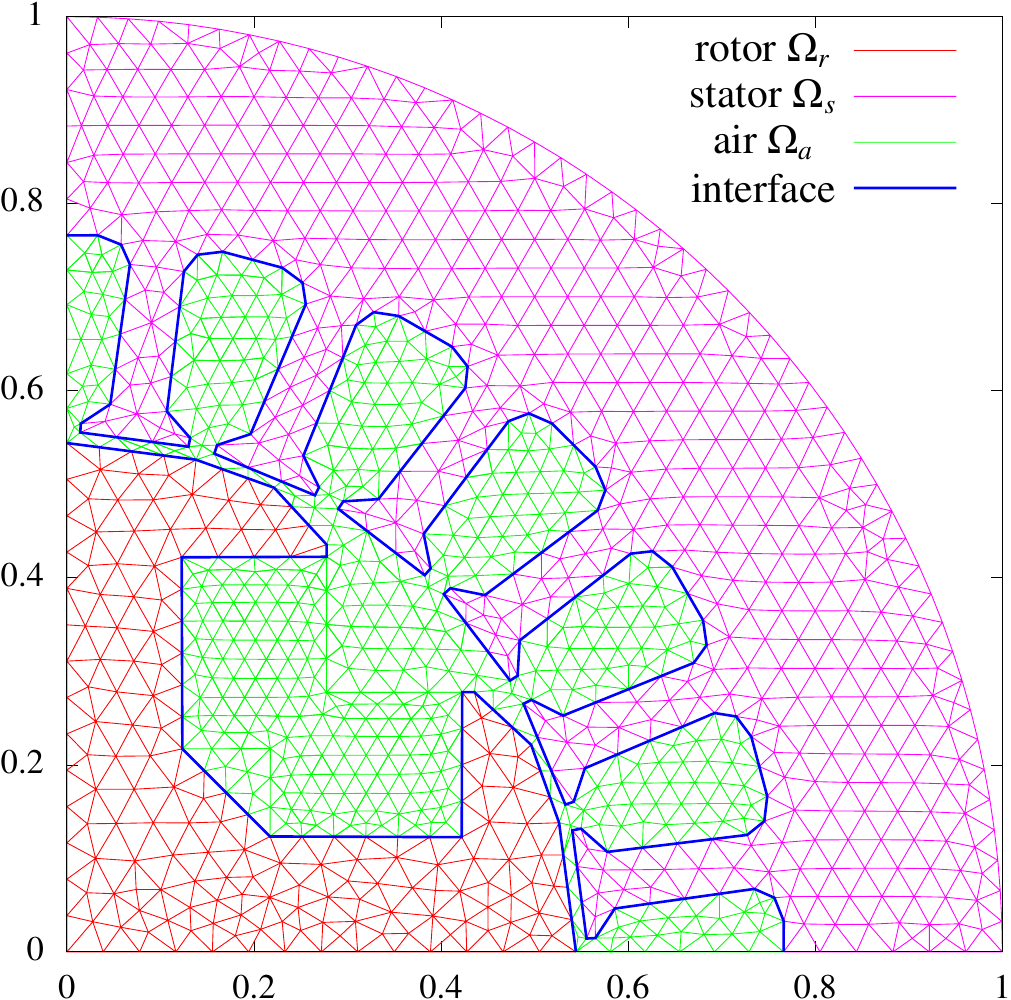}
    \hspace{0.01\textwidth}
    \includegraphics[height=0.34\textwidth]{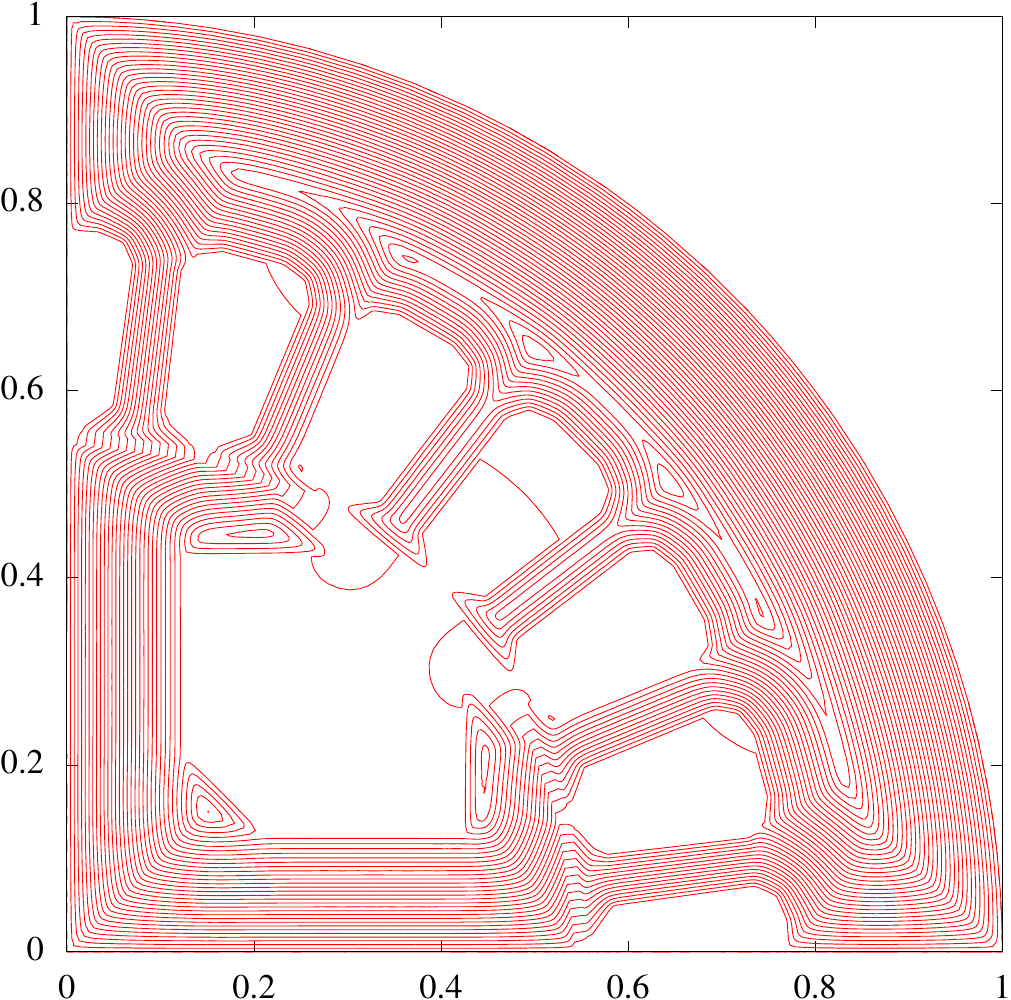}
    \hspace{0.01\textwidth}
    \includegraphics[height=0.34\textwidth]{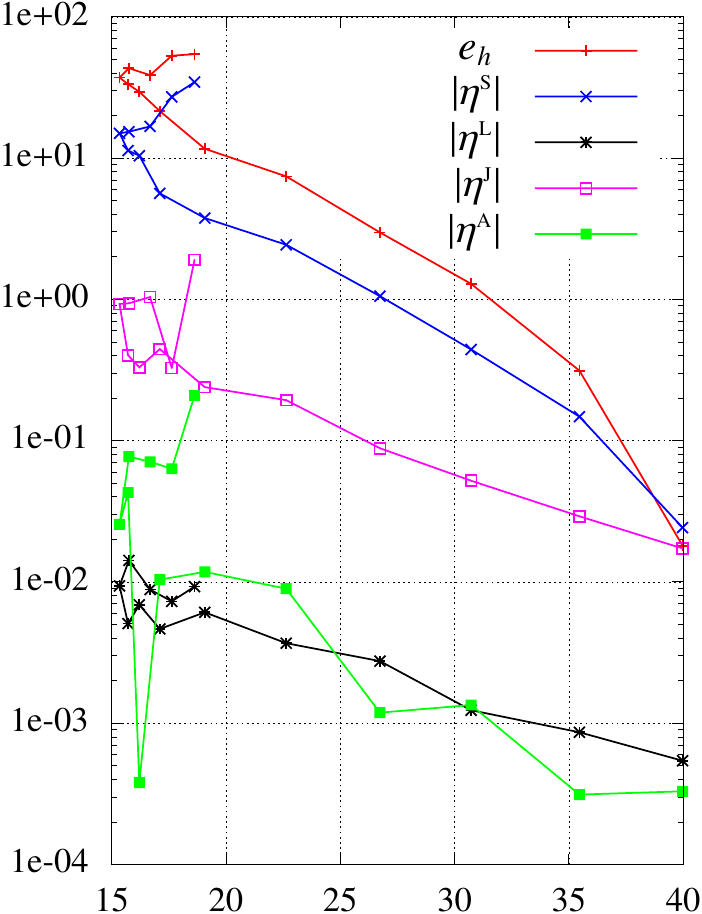}
  \end{center}
  \caption{Alternator, the computational domain with its components together with
    the initial mesh (left), the isolines of the primal solution (center) and the
    convergence of $e_h=|J(\u)-J(\uh)|$ and its estimators $|\etaS|$, $|\etaL|$, $|\etaJ|$
    and $|\etaA|$ with respect to $\DoF^{1/3}$  (right).}
  \label{fig:alter_geom}
\end{figure}

The last example follows from \cite{Glowinski74} where
the magnetic state in the cross-section of an alternator was solved numerically.
Due to the symmetry, only one quarter of the alternator is taken as the computational
domain $\Om:=\Om_s\cup\Om_s\cup\Om_a$;
see Figure~\ref{fig:alter_geom}, left, where the geometry of the domain is shown.
The alternator consists of the stator ($\Omega_s$) and rotor ($\Omega_r$)
with a gap filled by air ($\Omega_a$).

The problem is described by the Maxwell equations for the stationary magnetic field in the
form
\begin{subequations}
  \label{alt1}
  \begin{align}
    \label{alt1a}
    \rot H &= f \qquad \mbox{in }\Om, \\
    \label{alt1b}
    \Div B &= 0 \qquad \mbox{in }\Om, 
  \end{align}
\end{subequations}
where $H=(H_1, H_2)$ is the magnetic intensity field,
$B=(B_1, B_2)$ is the magnetic induction field and $f$ is the current density
(its component perpendicular to the plane of the computational domain).
The differential operators appearing in \eqref{alt1} are given by
$\rot H = (\pd H_2/\pd x_1,\, \pd H_1/\pd x_2)$ and
$\Div B = \nabla\cdot B = \pd B_1/\pd x_1 + \pd B_2/\pd x_2$ in two space dimensions.

\begin{table}
  \refA{\footnotesize
       \setlength{\tabcolsep}{2.5pt}
 \begin{tabular}{|cr|ccc|r|}
 \hline
 \multicolumn{6}{|c|}{$hp$-AMA            } \\
 \hline
  $\ell$ & $\DoF$ & 
  $ |e_h |$  &  $|\etaI|$ & $\ieff$ &time \\ 
 \hline
     0&  6462&  5.45E+01&  3.46E+01& 0.63& 313.5\\
     1&  5487&  5.31E+01&  2.71E+01& 0.51& 401.6\\
     2&  4648&  3.87E+01&  1.68E+01& 0.43& 501.4\\
     3&  3913&  4.33E+01&  1.54E+01& 0.35& 567.6\\
     4&  3617&  3.72E+01&  1.50E+01& 0.40& 603.3\\
 \hline
 \end{tabular}
 \hspace{1mm}
 \setlength{\tabcolsep}{2.5pt}
 \begin{tabular}{|cr|ccc|r|}
 \hline
 \multicolumn{6}{|c|}{$hp$-AMA            } \\
 \hline
  $\ell$ & $\DoF$ & 
  $ |e_h |$  &  $|\etaI|$ & $\ieff$ &time \\ 
 \hline
     5&  3883&  3.34E+01&  1.13E+01& 0.34& 678.7\\
     6&  4261&  2.94E+01&  1.04E+01& 0.35& 737.0\\
     7&  5010&  2.15E+01&  5.63E+00& 0.26& 819.0\\
     8&  6942&  1.16E+01&  3.77E+00& 0.32& 960.6\\
     9& 11597&  7.43E+00&  2.44E+00& 0.33&1224.0\\
 \hline
 \end{tabular}
 \hspace{1mm}
 \setlength{\tabcolsep}{2.5pt}
 \begin{tabular}{|cr|ccc|r|}
 \hline
 \multicolumn{6}{|c|}{$hp$-AMA            } \\
 \hline
  $\ell$ & $\DoF$ & 
  $ |e_h |$  &  $|\etaI|$ & $\ieff$ &time \\ 
 \hline
    10& 19107&  2.97E+00&  1.06E+00& 0.35&1717.4\\
    11& 29004&  1.29E+00&  4.42E-01& 0.34&2604.6\\
    12& 44588&  3.13E-01&  1.48E-01& 0.47&4210.5\\
    13& 63867&  1.81E-02&  2.43E-02& 1.34&8961.9\\
 \hline
 \end{tabular}

  }
  \caption{Alternator, error $e_h=|J(u)-J(\uh)|$,
    its estimate $|\etaI|$, effectivity index $\ieff$ and the computational time in seconds
    obtained by {\AMAhp} adaptive method for adaptive levels $\ell=0,\dots 13$.}
  \label{tab:alter}
\end{table}

Moreover, we consider the constitutive relation
\begin{align}
  \label{alt2}
  H(x) = \nu (x, |B(x)|^2) B(x), \quad x\in\Om,
\end{align}
where
\begin{align}
  \label{alt3}
  \nu(x, r) = 
  \begin{cases}
    \frac{1}{\mu_0} & \mbox{ for } x\in \Om_a, \\
    \frac{1}{\mu_0}\left(\alpha + (1-\alpha)\frac{r^4}{\beta+ r^4}\right) &
    \mbox{ for } x\in \Om_s\cup\Om_r. \\
  \end{cases}
\end{align}
The symbol
$\mu_0=1.256\times10^{-6}\, \mathrm{kg}\cdot\mathrm{m}\cdot\mathrm{A}^{-2}\cdot\mathrm{s}^{-2}$
denotes the permeability of
the vacuum and the material coefficients are $\alpha = 0.0003$, $\beta = 16000$ according
to \cite{Glowinski74}.
We consider the constant current density $f=5\times10^{4}\, \mathrm{A}\cdot\mathrm{m}^{-2}$.

Assuming that there exists a potential $u:\Om\to\R$ such that
$  B = \curl u = (\pd u/\pd x_2, - \pd u / \pd x_1)$,
equation \eqref{alt1b} is satisfied directly.
Obviously, $|B| = |\nabla u|$ and, therefore, \eqref{alt1a} together
with \eqref{alt2} gives
\begin{align}
  \label{alt5}
  f =  \rot H = \rot \nu (x, |B(x)|^2) B(x) =  \rot \nu (x, |\nabla u(x)|^2) \curl u (x)
  = -\nabla\cdot \left(\nu (x, |\nabla u(x)|^2) \nabla u\right).
\end{align}

Consequently, we have the following problem. Find $u:\Om\to \R$ such that
\begin{align}
  \label{alt6}
  -\nabla\cdot \left(\nu (x, |\nabla u(x)|^2) \nabla u\right) &= f\qquad \mbox{ in }\Om.
\end{align}
The homogeneous Dirichlet boundary condition is prescribed on $\Gamma$ for simplicity as in
\cite{Glowinski74}, but other options are possible.
We are interested in the total magnetic energy; hence, the target quantity is given by
\begin{align}
  \label{alt7}
  E = \frac 12 \int_{\Om} H(x) \cdot B(x)\dx = \frac 12\int_{\Om}  \nu (x, |B(x)|^2) |B(x)|^2\dx
  = \frac 12 \int_{\Om}  \nu (x, |\nabla u(x)|^2) |\nabla u(x)|^2\dx=:J(u).
\end{align}
Figure~\ref{fig:alter_geom}, center, shows the corresponding isolines of the primal solution
obtained on a fine grid where we obtained the reference value
$J(u) = 2\,664\, \mathrm{kg}\cdot\mathrm{m}^2\cdot\mathrm{s}^{-2}$.

\begin{figure} [t]
  \begin{center}
    \includegraphics[width=0.40\textwidth]{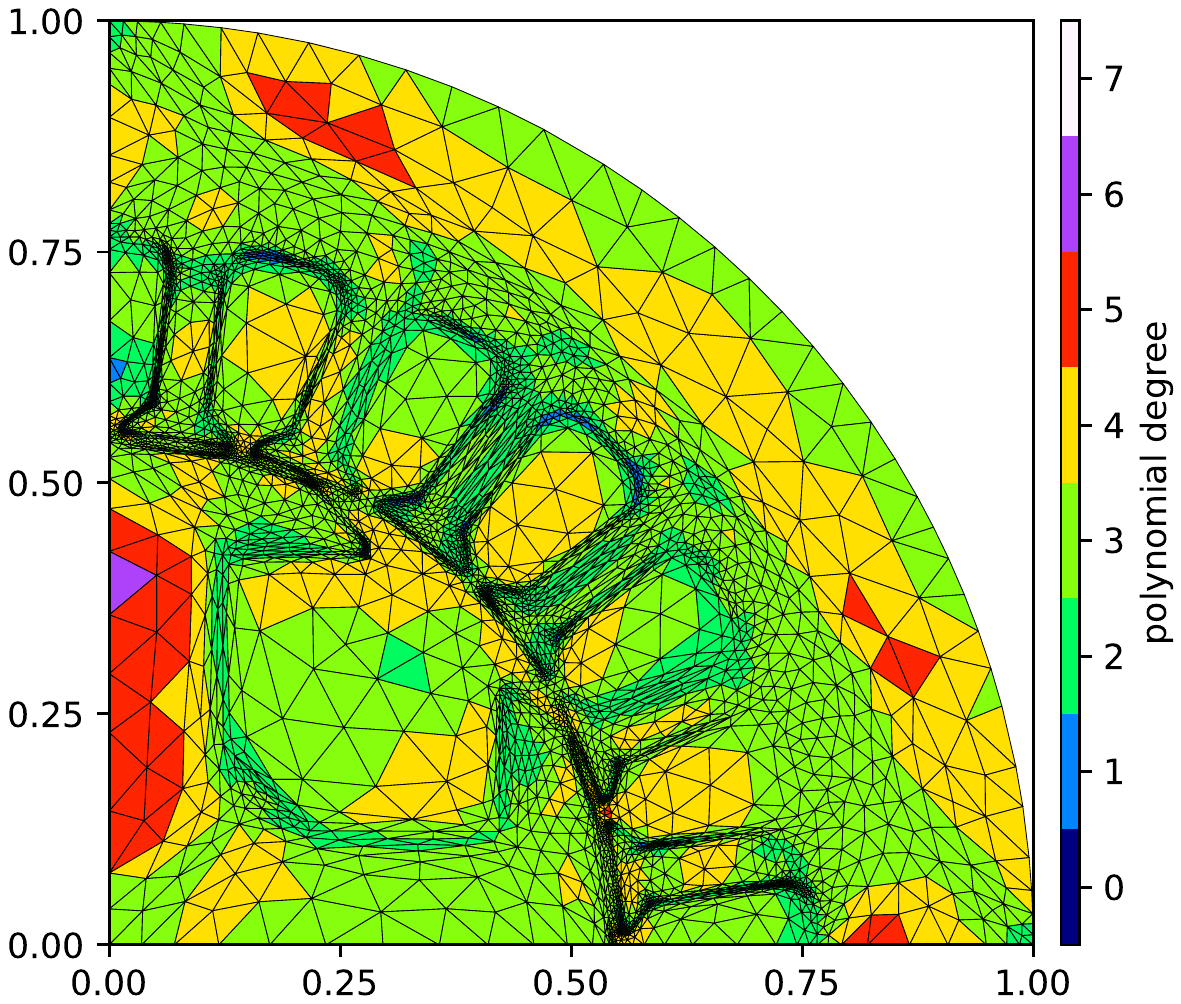}
    \hspace{0.06\textwidth}
    \includegraphics[width=0.40\textwidth]{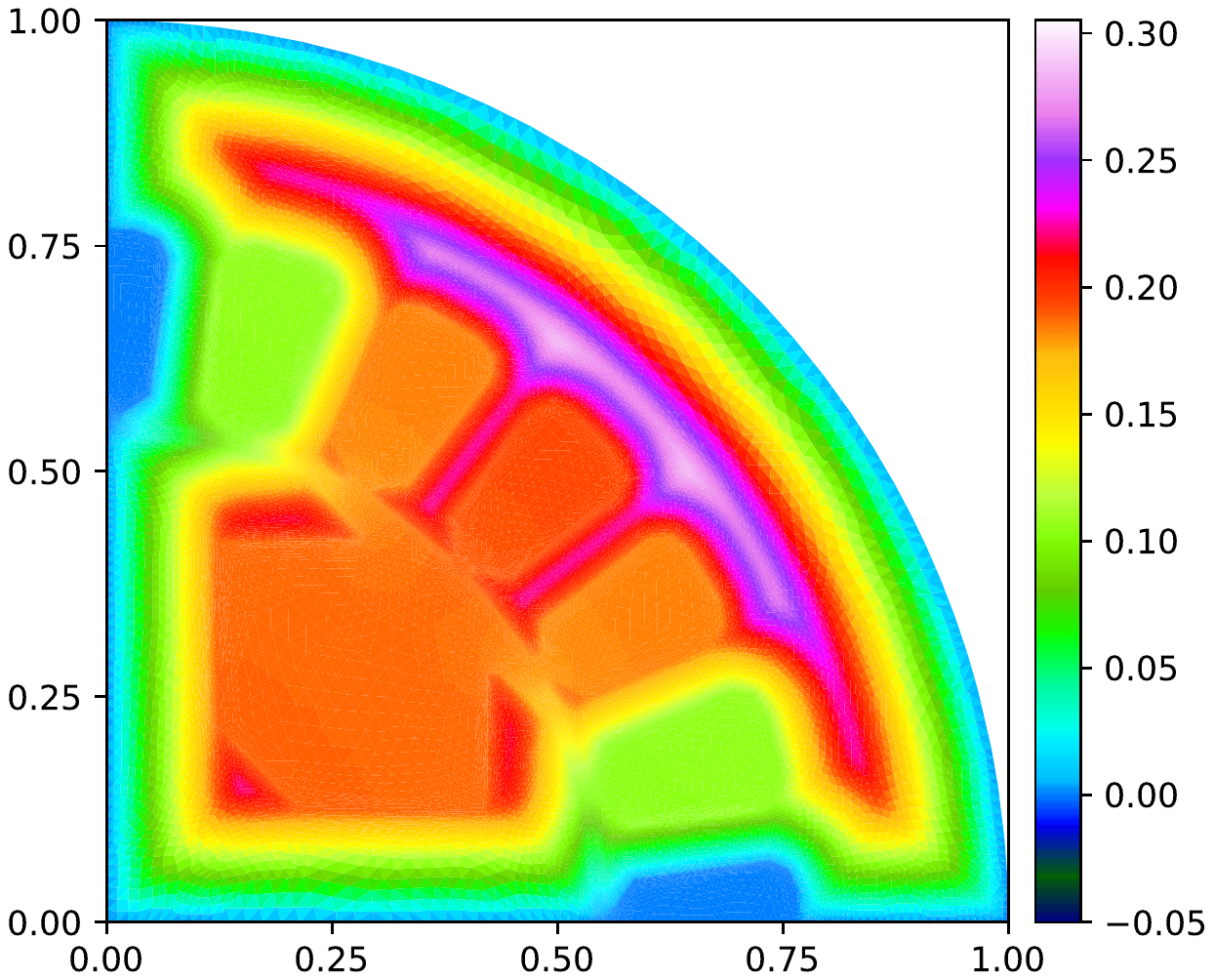}
  \end{center}
  \caption{Alternator, the final $hp$-mesh (left) and the corresponding magnetic potential (right).}
  \label{fig:alter_sol}
\end{figure}

We solve this problem with the {\AMAhp} technique.
Since function $\nu$ given by \eqref{alt3} differs for several
orders, the mesh adaptation technique must maintain the material interfaces.
Figure~\ref{fig:alter_geom}, right, shows the convergence of the error $e_h=|J(u)-\J(\uh)|$
and the error estimators $|\etaL|$, $|\etaJ|$, $|\etaS|$ and $|\etaA|$
with respect to $\DoF^{1/3}$, cf.~\eqref{CE10a}.
\refA{The values $e_h$, $\etaI$, effectivity index $\ieff$ and computational time
are shown in Table~\ref{tab:alter}.}
Again, we observe a reasonable approximation of the error
and the exponential rate of the convergence.
The final $hp$-grid and primal solution are plotted in Figure~\ref{fig:alter_sol}. A strong
refinement along the material interfaces is obvious.

\section{Conclusion}
\label{sec:concl}

We presented the framework of the goal-oriented error estimates for nonlinear problems
where the adjoint solution is based on the linearization of the primal weak formulation
used in the iterative solution of the corresponding algebraic systems.
We derived abstract error estimates consisting of three ingredients:
dual weighted residual, algebraic error and error arising from the linearization.
Then, employing a higher-order reconstruction, we proposed computable error estimates
and an adaptive algorithm for the numerical solution of nonlinear PDEs.
The presented numerical experiments demonstrate a reasonable approximation of the error
of the quantity of interest and also an exponential rate of convergence of the $hp$-adaptive method.
\refA{We are aware that the effectivity indexes are not enough close to the desired value
  around 1 but this is caused by an insufficient accuracy of the higher order approximation
  of the exact solutions of primal and dual problems on $hp$-anisotropic meshes.
  Nevertheless, the presented examples demonstrate a benefit of the use of such meshes.
  An improvement of the higher-order reconstruction will be a subject of further research.}


\begin{thebibliography}{10}
\providecommand{\url}[1]{{#1}}
\providecommand{\urlprefix}{URL }
\expandafter\ifx\csname urlstyle\endcsname\relax
  \providecommand{\doi}[1]{DOI~\discretionary{}{}{}#1}\else
  \providecommand{\doi}{DOI~\discretionary{}{}{}\begingroup
  \urlstyle{rm}\Url}\fi

\bibitem{Babuska-Suri90}
Babu{\v s}ka, I., Suri, M.: The $p$- and $hp$- versions of the finite element
  method. {A}n overview.
\newblock Comput. Methods Appl. Mech. Engrg. \textbf{80}, 5--26 (1990)

\bibitem{BalanWoopenMay16}
Balan, A., Woopen, M., May, G.: Adjoint-based $hp$-adaptivity on anisotropic
  meshes for high-order compressible flow simulations.
\newblock Comput. Fluids \textbf{139}, 47 -- 67 (2016)

\bibitem{RannacherBook}
Bangerth, W., Rannacher, R.: {Adaptive Finite Element Methods for Differential
  Equations}.
\newblock {Lectures in Mathematics. ETH Z{\"u}rich}. {Birkh\"auser Verlag}
  (2003)

\bibitem{BarrettLiu_NM93}
Barrett, J., Liu, W.: Finite element error analysis of a quasi-{N}ewtonian flow
  obeying the {C}arreau or power law.
\newblock Numer. Math. \textbf{64}(1), 433--453 (1993)

\bibitem{BeckerALL_CAMWA22}
Becker, R., Brunner, M., Innerberger, M., Melenk, J.M., Praetorius, D.:
  Rate-optimal goal-oriented adaptive {FEM} for semilinear elliptic {PDE}s.
\newblock Comput. Math. Appl. \textbf{118}, 18--35 (2022)

\bibitem{BeckerRannacher01}
Becker, R., Rannacher, R.: An optimal control approach to a-posteriori error
  estimation in finite element methods.
\newblock Acta Numerica \textbf{10}, 1--102 (2001)

\bibitem{BerroneSuli_IMA08}
Berrone, S., Süli, E.: Two-sided a posteriori error bounds for incompressible
  quasi-{N}ewtonian flows.
\newblock IMA J. Numer. Anal. \textbf{28}(2), 382--421 (2008)

\bibitem{ChaillouSuri_CMAME06}
Chaillou, A., Suri, M.: Computable error estimators for the approximation of
  nonlinear problems by linearized models.
\newblock Comput. Methods Appl. Mech. Engrg. \textbf{196}(1-3), 210--224 (2006)

\bibitem{CongreveALL_IMA13}
Congreve, S., Houston, P., Süli, E., Wihler, T.: Discontinuous {G}alerkin
  finite element approximation of quasilinear elliptic boundary value problems
  {II}: {S}trongly monotone quasi-{N}ewtonian flows.
\newblock IMA J. Numer. Anal. \textbf{33}(4), 1386--1415 (2013)

\bibitem{CongreveWihler_JCAM17}
Congreve, S., Wihler, T.P.: Iterative {G}alerkin discretizations for strongly
  monotone problems.
\newblock J. Comput. Appl. Math. \textbf{311}, 457 -- 472 (2017)

\bibitem{Demko02}
Demkowicz, L., Rachowicz, W., Devloo, P.: A fully automatic $hp$-adaptivity.
\newblock J. Sci. Comput. \textbf{17}(1-4), 117--142 (2002)

\bibitem{DiStolfoALL_JNM19}
Di~Stolfo, P., Rademacher, A., Schröder, A.: Dual weighted residual error
  estimation for the finite cell method.
\newblock Journal of Numerical Mathematics \textbf{27}(2), 101--122 (2019)

\bibitem{GO_nonlinear}
Dolej{\v s}{\'\i}, V., Barto{\v s}, O., Roskovec, F.: Goal-oriented mesh
  adaptation method for nonlinear problems including algebraic errors.
\newblock Comput. Math. Appl. \textbf{93}, 178--198 (2021)

\bibitem{DGM-book}
Dolej{\v s}{\'\i}, V., Feistauer, M.: Discontinuous {G}alerkin Method --
  Analysis and Applications to Compressible Flow.
\newblock Springer Series in Computational Mathematics 48. Springer, Cham
  (2015)

\bibitem{AMA-book}
Dolej{\v s}{\'\i}, V., May, G.: Anisotropic $hp$-Mesh Adaptation Methods.
\newblock Birkh\"auser (2022)

\bibitem{ESCO-16}
Dolej{\v s}{\'\i}, V., May, G., Roskovec, F., Solin, P.: Anisotropic $hp$-mesh
  optimization technique based on the continuous mesh and error models.
\newblock Comput. Math. Appl. \textbf{74}, 45--63 (2017)

\bibitem{hp_reconstr}
Dolej{\v s}{\'\i}, V., Solin, P.: $hp$-discontinuous {G}alerkin method based on
  local higher order reconstruction.
\newblock Appl. Math. Comput. \textbf{279}, 219--235 (2016)

\bibitem{DolTich_BiCG}
Dolej{\v s}{\'\i}, V., Tich{\'y}, P.: On efficient numerical solution of linear
  algebraic systems arising in goal-oriented error estimates.
\newblock Journal of Scientific Computing \textbf{83}(5) (2020)

\bibitem{EndtmayerLangerWisk_SISC20}
Endtmayer, B., Langer, U., Wick, T.: Two-side a posteriori error estimates for
  the dual-weighted residual method.
\newblock SIAM Journal on Scientific Computing \textbf{42}(1), A371--A394
  (2020)

\bibitem{ErnVohralik_SISC13}
Ern, A., Vohral{\'\i}k, M.: Adaptive inexact {N}ewton methods with a posteriori
  stopping criteria for nonlinear diffusion {PDE}s.
\newblock SIAM J. Sci. Comput. \textbf{35}(4), A1761--A1791 (2013)

\bibitem{FidkowskiDarmofal_AIAA11}
Fidkowski, K., Darmofal, D.: Review of output-based error estimation and mesh
  adaptation in computational fluid dynamics.
\newblock AIAA Journal \textbf{49}(4), 673--694 (2011)

\bibitem{GileSuli02}
Giles, M., S{\"u}li, E.: Adjoint methods for {PDE}s: a posteriori error
  analysis and postprocessing by duality.
\newblock Acta Numerica \textbf{11}, 145--236 (2002)

\bibitem{Glowinski74}
Glowinski, R., Marrocco, A.: Analyse num\'erique du champ magnetique d’un
  alternateur par elements finis et sur-relaxation ponctuelle non lineaire.
\newblock Comput. Methods Appl. Mech. Engrg. \textbf{3}, 55--85 (1974)

\bibitem{HaberlAll_NM21}
Haberl, A., Praetorius, D., Schimanko, S., Vohralik, M.: Convergence and
  quasi-optimal cost of adaptive algorithms for nonlinear operators including
  iterative linearization and algebraic solver.
\newblock Numer. Math. \textbf{147}(3), 679--725 (2021)

\bibitem{Hartmann2007Adjoint}
Hartmann, R.: {Adjoint Consistency Analysis of Discontinuous {G}alerkin
  Discretizations}.
\newblock SIAM J. Numer. Anal. \textbf{45}(6), 2671--2696 (2007)

\bibitem{HH06:SIPG2}
Hartmann, R., Houston, P.: Symmetric interior penalty {DG} methods for the
  compressible {N}avier-{S}tokes equations {II}: {G}oal-oriented a posteriori
  error estimation.
\newblock Int. J. Numer. Anal. Model. \textbf{3}, 141--162 (2006)

\bibitem{HeidAll_CMAM21}
Heid, P., Praetorius, D., Wihler, T.P.: Energy contraction and optimal
  convergence of adaptive iterative linearized finite element methods.
\newblock Comput. Meth. Aappl. Math. \textbf{21}(2, SI), 407--422 (2021)

\bibitem{HeidWihler_CAL20}
Heid, P., Wihler, T.: On the convergence of adaptive iterative linearized
  {G}alerkin methods.
\newblock Calcolo \textbf{57}(3) (2020)

\bibitem{houston-robson-suli}
Houston, P., Robson, J., S{\"u}li, E.: Discontinuous {G}alerkin finite element
  approximation of quasilinear elliptic boundary value problems {I}: {T}he
  scalar case.
\newblock IMA J. Numer. Anal. \textbf{25}, 726--749 (2005)

\bibitem{ListRadu_ComputGeosc16}
List, F., Radu, F.A.: {A study on iterative methods for solving {R}ichards'
  equation}.
\newblock Comput. Geosci. \textbf{{20}}({2}), {341--353} ({2016})

\bibitem{LoseilleAlauzet11a}
Loseille, A., Alauzet, F.: Continuous mesh framework part {I}: well-posed
  continuous interpolation error.
\newblock SIAM J. Numer. Anal. \textbf{49}(1), 38--60 (2011)

\bibitem{LoseilleDervieuxAlauzet_JCP10}
Loseille, A., Dervieux, A., Alauzet, F.: Fully anisotropic goal-oriented mesh
  adaptation for 3{D} steady {E}uler equations.
\newblock J. Comput. Phys. \textbf{229}(8), 2866--2897 (2010)

\bibitem{MakridakisNochetto2003}
Makridakis, C., Nochetto, R.H.: {Elliptic reconstruction and a posteriori error
  estimates for parabolic problems}.
\newblock {SIAM J. Numer. Anal.} \textbf{41}(4), 1585--1594 (2003)

\bibitem{MallikVohralikYousef_JCAM20}
Mallik, G., Vohral{\'\i}k, M., Yousef, S.: Goal-oriented a posteriori error
  estimation for conforming and nonconforming approximations with inexact
  solvers.
\newblock Journal of Computational and Applied Mathematics \textbf{366} (2020)

\bibitem{RaduALL_JCAM15}
Radu, F.A., Nordbotten, J.M., Pop, I.S., Kumar, K.: {A robust linearization
  scheme for finite volume based discretizations for simulation of two-phase
  flow in porous media}.
\newblock {J. Comput. Appl. Math.} \textbf{{289}}, {134--141} ({2015})

\bibitem{Rannacher201323}
Rannacher, R., Vihharev, J.: Adaptive finite element analysis of nonlinear
  problems: Balancing of discretization and iteration errors.
\newblock J. Numer. Math. \textbf{21}(1), 23--61 (2013)

\bibitem{RichterWick_JCAM15}
Richter, T., Wick, T.: Variational localizations of the dual weighted residual
  estimator.
\newblock J. Comput. Appl. Math. \textbf{279}, 192 -- 208 (2015)

\bibitem{Schwab-book}
Schwab, C.: {$p$}- and {$hp$}-Finite Element Methods.
\newblock Clarendon Press, Oxford (1998)

\bibitem{SolDemko04}
{\v S}ol{\'\i}n, P., Demkowicz, L.: Goal-oriented $hp$-adaptivity for elliptic
  problems.
\newblock Comput. Methods Appl. Mech. Engrg. \textbf{193}, 449--468 (2004)

\bibitem{ZeidlerIIB}
Zeidler, E.: Nonlinear functional analysis and its applications. {II}/{B},
  {N}onlinear monotone operators.
\newblock New York, Springer (1985)

\end{thebibliography}

\end{document}